\documentclass[12pt]{article}
\usepackage{amsmath,amsthm,amsfonts,amssymb, amscd,amssymb,latexsym}
\usepackage{graphicx}
\usepackage{hyperref}
\usepackage{srcltx}
\usepackage[usenames]{color}
\usepackage[matrix,arrow,curve]{xy}

\newtheorem{theorem}{Theorem}[section]
\newtheorem{lemma}[theorem]{Lemma}
\theoremstyle{definition}

\newtheorem{corollary}[theorem]{Corollary}

\newtheorem{definition}[theorem]{Definition}
\newtheorem{example}[theorem]{Example}

\newtheorem{remark}[theorem]{Remark}
\newtheorem{proposition}[theorem]{Proposition}

\renewcommand{\L}{{\mathtt{L}}}
\def\C{{\mathcal C}}
\def\A{{\mathcal{A}}}
\def\B{{\mathcal{B}}}
\def\Ce{{\mathcal{C}}}
\def\M{{\mathcal{M}}}

\def\E{{\mathcal{E}}}
\def\T{{\mathcal{T}}}

\newcommand{\Th}{{\mathrm{Th}}}
\newcommand{\Tr}{{\mathrm{T}}}
\newcommand{\At}{{\mathrm{At}}}
\newcommand{\Mod}{{\mathrm{Mod}}}

\newcommand{\uw}{{\mathrm{u}_\omega}}
\newcommand{\qw}{{\mathrm{q}_\omega}}
\newcommand{\us}{{\mathrm{u}_S}}
\newcommand{\qs}{{\mathrm{q}_S}}
\newcommand{\ur}{{\mathrm{u}_R}}

\newcommand{\K}{{\mathbf{K}}}

\renewcommand{\S}{{\mathbf{S}}}
\renewcommand{\P}{{\mathbf{P}}}
\newcommand{\Pw}{{\mathbf{P_{\!\! \omega}}}}
\newcommand{\Ps}{{\mathbf{P_{\!\! s}}}}
\newcommand{\Pf}{{\mathbf{P_{\! f}}}}
\newcommand{\Pu}{{\mathbf{P_{\!\! u}}}}

\newcommand{\Ld}{{\underrightarrow{\mathbf{L}}}}
\newcommand{\Ls}{{\underrightarrow{\mathbf{L}}_{\! s}}}

\newcommand{\Lo}{{\mathbf{L_{fg}}}}
\newcommand{\e}{{\mathbf{\,_e}}}

\newcommand{\UC}{{\mathbf{U}}}
\newcommand{\QC}{{\mathbf{Q}}}
\newcommand{\EN}{{\mathbf{N}}}

\newcommand{\pvar}{{\mathbf{Pvar}}}

\newcommand{\qvar}{{\mathbf{Qvar}}}
\newcommand{\ucl}{{\mathbf{Ucl}}}
\newcommand{\Res}{{\mathbf{Res}}}
\newcommand{\Dis}{{\mathbf{Dis}}}

\newcommand{\Rad}{{\mathrm{Rad}}}
\newcommand{\V}{{\mathrm{V}}}

\title{Algebraic geometry over algebraic structures III: Equationally Noetherian property and compactness}

\author{E.\,Daniyarova, A.\,Myasnikov, V.\,Remeslennikov}

\begin{document}
\maketitle

\begin{abstract}
In this paper we discuss some special generalizations of
equationally Noetherian property which naturally arise in the
universal algebraic geometry. We introduce weakly equationally
Noetherian, $\qw$-compact, $\uw$-compact, and weakly $\uw$-compact
algebras and then examine properties of such algebras. Also we
consider the connections between five classes: the class of
equationally Noetherian algebras, the class of weakly equationally
Noetherian algebras, the class of $\uw$-compact algebras, the
class of weakly $\uw$-compact algebras, and the class of
$\qw$-compact algebras.


\end{abstract}

\tableofcontents

\section{Introduction}

This paper deals with the {\em universal algebraic geometry}. The
universal algebraic geometry is a young branch of mathematics. The
subject of universal algebraic geometry lies in the solutions of
systems of equations over an arbitrary algebraic structure.

Investigations in universal algebraic geometry were started in
works by B.\,I.\,Plotkin~\cite{Plot1, Plot2, Plot3} and papers on
algebraic geometry over groups by G.\,Baumslag,
O.\,G.\,Kharlampovich, A.\,G.\,Myasnikov, and
V.\,N.\,Remeslennikov~\cite{KM1, KM2, KM3, KM4}. After that there
were a lot of papers on algebraic geometry over concrete groups,
algebras, monoids and so on. Among them there are the famous works
by O.\,G.\,Kharlampovich, A.\,G.\,Myasnikov~\cite{KM1, KM2, KM3,
KM4} and Z.\,Sela~\cite{Sela1, Sela2, Sela3} on algebraic geometry
over free groups.

In recent years we have achieved more general and systematic point
of view on the universal algebraic geometry as on a formalized
theory. In this respect we have started a series of works on
universal algebraic geometry. This paper is the third one of that
series along with~\cite{DMR1, DMR2}.

According to~\cite{Gorbunov, Malcev, Marker}, in~\cite{DMR1} we
give a framework of universal algebra and model theory as much as
we need it in universal algebraic geometry. At the same time we
discuss how notions and ideas from model theory work in universal
algebraic geometry. In~\cite{DMR2} we introduce the foundation of
universal algebraic geometry, basic definitions and constructions
of the algebraic geometry over an arbitrary algebraic structure
$\B$.

This paper is supposed to be read after the previous
ones~\cite{DMR1, DMR2}, however for the sake of convenience we
present in here some of the most essential notations and
definitions (see Section~\ref{sec:preliminaries}).

We consider only first-order functional languages (signatures).
Recall that algebraic structures in a functional language are
called {\em algebras}. Typically we denote algebraic structures by
capital calligraphic letters ($\A, \B, \Ce, \ldots$) and their
universes by the corresponding capital Latin letters ($A, B, C,
\ldots$).

The main results of papers~\cite{DMR1, DMR2} are so-called the
Unification Theorems (Theorem~A and Theorem~C) which give a
description of coordinate algebras by means of several languages.

\noindent {\bf Theorem~A.} {\it Let $\B$ be an equationally
Noetherian algebra in a functional language $\L$. Then for a
finitely generated algebra $\Ce$ of $\L$ the following conditions
are equivalent:
\begin{enumerate}
\item [1)] $\Th_{\forall} (\B) \subseteq \Th_{\forall} (\Ce)$, i.e., $\Ce \in \ucl(\B)$;
\item [2)] $\Th_{\exists} (\B) \supseteq \Th_{\exists} (\Ce)$;
\item [3)] $\Ce$ embeds into an ultrapower of $\B$;
\item [4)] $\Ce$ is discriminated by $\B$;
\item [5)] $\Ce$ is a limit algebra over $\B$;
\item [6)] $\Ce$ is an algebra defined by a complete atomic type
in the theory $\Th _{\forall} (\B)$ in $\L$;
\item [7)] $\Ce$ is the coordinate algebra of an
irreducible algebraic set over $\B$ defined by a system of
equations in the language $\L$.
\end{enumerate}
}

\noindent {\bf Theorem~C.} {\it Let $\B$ be an equationally
Noetherian algebra in a functional language $\L$. Then for a
finitely generated algebra $\Ce$ of $\L$ the following conditions
are equivalent:
\begin{enumerate}
\item [1)] $\Ce \in \qvar(\B)$, i.e., $\Th_{\rm qi} (\B) \subseteq \Th_{\rm qi} (\Ce)$;
\item [2)] $\Ce \in \pvar(\B)$;
\item [3)] $\Ce$ embeds into a direct power of $\B$;
\item [4)] $\Ce$ is separated by $\B$;
\item [5)] $\Ce$ is a subdirect product of a finitely many limit algebras over $\B$;
\item [6)] $\Ce$ is an algebra defined by a complete atomic type in the theory $\Th _{\rm qi} (\B)$ in $\L$;
\item [7)] $\Ce$ is the coordinate algebra of an algebraic set over $\B$ defined by a system of equations in the language $\L$.
\end{enumerate}
}

Note that items 5) in both Theorem~A and Theorem~C give a
description of coordinate algebras by means of limit algebras. The
limit algebraic structures (groups, as the rule) become the object
of intense interest in modern algebra~\cite{CG, GS, Groves1,
Groves1-2, Groves2-1, Groves2, G1}.

\bigskip

Theorems~A and~C are formulated for so-called equationally
Noetherian algebras (the definition see in
Section~\ref{sec:preliminaries}). Equationally Noetherian algebras
possess the best opportunity to study the algebraic geometry over
them. If a given algebra $\B$ is equationally Noetherian then we
have an advantage when investigating the algebraic geometry over
$\B$. In this case we may use:
\begin{itemize}
\item[(i)] Unification Theorems;
\item[(ii)] the decomposition of any algebraic set over $\B$ into a
finite union of irreducible algebraic sets (Theorem~\ref{irr}
below);
\item[(iii)] the possibility to study only finite system
of equations;
\item[(iv)] and some more results~\cite{DMR1, DMR2}.
\end{itemize}

In the case when a given algebra $\B$ is not equationally
Noetherian we lose some results for equationally Noetherian
algebras, while some of them may remain in force. In this paper we
introduce four generalizations of the equationally Noetherian
property which naturally arise in universal algebraic geometry.
These are
\begin{itemize}
\item[($\EN'$):] weak equationally Noetherian property that retains~(iii);
\item[($\QC\,$):] $\qw$-compactness that retains Unification Theorem~C;
\item[($\UC\,$):] $\uw$-compactness that retains Unification Theorems~A and~C;
\item[($\UC'$):] and weak $\uw$-compactness that retains~(iv), namely, some weak form of Unification Theorem~A.
\end{itemize}

We denote by $\EN$ the class of all equationally Noetherian
algebras in a given functional language $\L$. By $\EN'$, $\QC$,
$\UC$, $\UC'$, correspondingly, we denote the classes of algebras
with properties above. The picture of connections between classes
$\QC$, $\UC$, $\UC'$, $\EN'$ and $\EN$ is presented in
Section~\ref{sec:picture}.

\medskip

There exist several equivalent approaches to $\qw$- and
$\uw$-compact algebras. We introduce them in
Section~\ref{sec:compact}. One of these approaches rises from some
ideas of model theory. It relates to the Compactness Theorem and
the notion of compact algebra.

Recall that a set of formulas $T$ in a language $\L$ is called
{\em satisfiable} in a class $\K$ of algebraic structures in $\L$
(or $T$ is {\em realized} in $\K$) if one can assign some elements
from a particular algebraic structure from $\K$ as values to the
variables which occur in $T$ in such a way that all formulas from
$T$ become true. The set $T$ is called {\em finitely satisfiable}
in $\K$ if every finite subset of $T$ is realized in $\K$.

\medskip

\noindent {\bf Compactness Theorem} (K.\,F.\,G\"{o}del,
A.\,I.\,Malcev~\cite{Gorbunov}){\bf.} {\it If a set of first-order
formulas $T$ in a language $\L$ is finitely satisfiable in a class
$\K$ of algebraic structures in $\L$, then $T$ is satisfiable in
an ultraproduct of structures from $\K$.}

Class $\K$ is called {\em compact} if every finitely satisfiable
in $\K$ set of formulas $T$ is satisfiable in $\K$. This
definition occurs in the book by V.\,A.\,Gorbunov~\cite{Gorbunov}.
It is natural to name an algebraic structure $\B$ compact if the
class $\{\B\}$ is compact. However, according to
W.\,Hodges~\cite{Hodges}, algebraic $\L$-structure is called
compact if its universe is a Hausdorff topological space, in such
a way that each function from $\L$ is interpreted by a continuous
function. The same algebraic structures appear in~\cite{Gorbunov}
under the name of topologically compact structures.

Trying to avoid an ambiguity we call an $\L$-algebra $\B$ {\em
logically compact} if every finitely satisfiable in $\B$ set of
formulas $T$ in the language $\L$ is satisfiable in $\B$. When we
modify this definition and consider only special types of sets of
formulas $T$ we get definitions of special compactness, such as
$\qw$- and $\uw$-compactness. Short review of the history of
``$\qw$-compact'' notion is represented in
Subsection~\ref{subsec:compact}.

\bigskip

First and foremost in this article we generalize the Unification
Theorems to $\uw$- and $\qw$-compact algebras. In
Section~\ref{sec:unification_theorems_compact} we give geometric
definitions of $\uw$- and $\qw$-compactness. In
Subsection~\ref{subsec:A} we prove that Theorem~A is true for any
$\uw$-compact algebra $\B$ and every algebra $\B$ which satisfies
Theorem~A is $\uw$-compact. The similar result that connects
$\qw$-compact algebras and Theorem~C is presented in
Subsection~\ref{subsec:C}. In Subsection~\ref{subsec:weak} for
weakly $\uw$-compact algebras we formulate and prove a weak analog
of Theorem~A.

Section~\ref{sec:compact} is devoted to $\qw$- and $\uw$-compact
algebras. In Subsection~\ref{subsec:compact} we put definitions of
$\qw$- and $\uw$-compact algebras in different equivalent forms
and prove the equivalence of them in
Subsection~\ref{subsec:proof}. For $\uw$-compact algebra $\B$
Unification Theorems give a global view to all (irreducible)
coordinate algebras over $\B$. However, it may happen that one has
no $\uw$-compact property but some ``local $\uw$-compact
property'' which gives result of Theorem~A for a certain algebra
$\Ce$ (not for all $\Ce$). This idea is developed in
Subsection~\ref{subsec:local}.

In Section~\ref{sec:weak} we discuss weak properties: weak
equationally Noetherian property (Subsection~\ref{subsec:weakEN})
and weak $\uw$-compactness (Subsection~\ref{subsec:weakcompact}).
In Subsection~\ref{subsec:logirr} we introduce logically
irreducible algebraic sets. Those sets naturally arise as
generalization of irreducible ones. In particularly, we show that
the notions of irreducible algebraic set and logically irreducible
algebraic set over an algebra $\B$ coincide if and only if $\B$ is
weakly $\uw$-compact.

In the last Section~\ref{sec:extension} we continue discussion
about connections between $\uw$- and $\qw$-compact algebras with
the Compactness Theorem and corresponding technique from the model
theory. By the way, we construct $\uw$-compact elementary
extension for an arbitrary algebra $\B$.

\section{Preliminaries}\label{sec:preliminaries}

In this section we remind basic notions and facts from universal
algebraic geometry according to~\cite{DMR1, DMR2}.

Let $\L$ be a first-order functional language, $X = \{x_1, x_2,
\ldots, x_n\}$ a finite set of variables, $\Tr_\L(X)$ the set of
all terms of $\L$ with variables in $X$, $\T_\L(X)$ the absolutely
free $\L$-algebra with basis $X$ and $\At_{\L}(X)$ the set of all
atomic formulas of $\L$ with variables in $X$.

In universal algebraic geometry atomic formulas from $\At_{\L}(X)$
are named {\em equations} in $\L$ and subsets $S \subseteq
\At_{\L}(X)$ are named {\em systems of equations} in the language
$\L$.

For a system of equations $S \subseteq \At_{\L}(X)$ and an algebra
$\B$ in the language $\L$ we denote by $\V _\B(S)$ the set of all
solutions of $S$ in $\B$:
$$
\V _\B(S)=\{(b_1,\ldots,b_n)\in B^n \; \vert\quad \B\models
(t(b_1,\ldots,b_n)=s(b_1,\ldots,b_n))\quad \forall \;(t=s)\in S\}.
$$
It is called the {\em algebraic set} over $\B$ defined by the
system $S$. If $S$ contains of only one equation $(t=s)$ we write
$\V_\B(t=s)$ instead of $\V_\B(\{(t=s)\})$.

Algebraic set is {\em irreducible} if it is not a finite union of
proper algebraic subsets; otherwise it is {\em reducible}. The
empty set is not considered to be irreducible. Hence, according to
R.\,Hartshorne~\cite{Hartshorne}, all irreducible algebraic sets
are non-empty in our paper.

Two systems $S_1, S_2 \subseteq \At_{\L}(X)$ are {\em equivalent}
over $\B$ if $\V_\B(S_1)=\V_\B(S_2)$. The {\em radical}
$\Rad_\B(S)$ of a system of equations $S \subseteq \At_{\L}(X)$ is
the maximal system which is equivalent to $S$ over $\B$. It is
also called the radical of algebraic set $Y=\V_\B(S)$ and denoted
by $\Rad(Y)$. By $[S]$ we denote the congruent closure of $S$,
i.e., the least congruent subset of $\At_{\L}(X)$ that contains
$S$.

By $\Phi _{\mathrm{qf},\L} (X)$ we denote the set of all
quantifier-free formulas in $\L$ with variables in $X$. We say
that a formula $\phi \in \Phi _{\mathrm{qf},\L} (X)$ is a {\em
consequence} of a system of equations $S\subseteq \At_\L (X)$ over
an $\L$-algebra $\B$, if $\B\models \phi (b_1,\ldots,b_n)$ for all
$(b_1,\ldots,b_n)\in \V_\B (S)$. For example, an atomic formula
$(t=s)$, $t,s\in \Tr_\L(X)$, is a consequence of $S$ over $\B$ if
and only if $(t=s)\in \Rad_\B(S)$.

For an arbitrary algebraic set $Y \subseteq B^n$ over $\B$ the
radical $\Rad(Y)$ defines the congruence $\theta_{\Rad(Y)}$ on
$\T_\L(X)$:
$$
t_1\sim_{\theta_{\Rad(Y)}}t_2 \quad \Longleftrightarrow \quad
(t_1=t_2)\in \Rad(Y), \quad t_1,t_2\in \Tr_\L(X).
$$
The factor-algebra $\Gamma (Y)=\T_\L(X)/\theta_{\Rad(Y)}$ is
called the {\em coordinate algebra} of the algebraic set $Y$.

Let $Y\subseteq B^n$ and $Z\subseteq B^m$ be algebraic sets over
$\B$. One has $\Gamma (Y)\cong\Gamma (Z)$ if and only if algebraic
sets $Y$ and $Z$ are isomorphic (we omit here the definition of
isomorphism between algebraic sets). Isomorphic algebraic sets are
irreducible and reducible simultaneously.

We say that an $\L$-algebra $\Ce$ is a {\em coordinate algebra}
over $\B$ if $\Ce \cong \Gamma (Y)$ for some algebraic set $Y$
over $\B$, and $\Ce$ is an {\em irreducible coordinate algebra}
over $\B$ if $\Ce \cong \Gamma (Y)$ for some irreducible algebraic
set $Y$ over $\B$.

One of the principal goals of algebraic geometry over a given
algebraic structure $\B$ is the problem of classification of
algebraic sets over $\B$ up to isomorphism. This problem is
equivalent to the problem of classification of coordinate algebras
of algebraic sets over $\B$. Also it is important to classify
coordinate algebras of irreducible algebraic sets over $\B$.
Formulated in Introduction Unification Theorems~A and~C are very
useful for solution of those problems.

In Theorems~A and~C we claim an algebra $\B$ is equationally
Noetherian. Thus, let us remind that an $\L$-algebra $\B$ is
called {\em equationally Noetherian}, if for every finite set $X$
and every system of equations $S \subseteq \At_{\L}(X)$ there
exists a finite subsystem $S_0\subseteq S$ such that
$\V_\B(S_0)=\V_\B(S)$. Properties of equationally Noetherian
algebras are discussed in~\cite{DMR1, DMR2}.

An $\L$-algebra $\C$ is {\it separated} by $\L$-algebra $\B$ if
for any pair of non-equal elements $c_1, c_2 \in C$ there is a
homomorphism $h \colon \C \to \B$ such that $h(c_1)\ne h(c_2)$. An
algebra $\C$ is {\em discriminated} by $\B$ if for any finite set
$W$ of elements from $C$ there is a homomorphism $h \colon\C \to
\B$ whose restriction onto $W$ is injective. We are interested in
a familiar form of results, so it is useful to put by definition
that the trivial algebra $\E$ is separated by an algebra $\B$
anyway, and $\E$ is discriminated by $\B$ if and only if $\B$ has
a trivial subalgebra.

The definitions of limit algebras and algebras defined by complete
atomic types need a large introduction, so we omit them
(see~\cite{DMR1}).

In this paper we use some operators which image a class $\K$ of
$\L$-algebras into another one. For the sake on convenience we
collect here the list of all these operators:\\
$\S (\K)$~--- the class of subalgebras of algebras from $\K$;\\
$\P (\K)$~--- the class of direct products of algebras from
$\K$;\\ $\Pw (\K)$~--- the class of finite direct products of
algebras from $\K$;\\
$\Ps (\K)$~--- the class of subdirect products of algebras from
$\K$;\\
$\Pf (\K)$~--- the class of filterproducts of algebras from
$\K$;\\
$\Pu (\K)$~--- the class of ultraproducts of algebras from $\K$;\\
$\Ld (\K)$~--- the class of direct limits of algebras from $\K$;\\
$\Ls (\K)$~--- the class of epimorphic direct limits of algebras
from $\K$;\\
$\Lo(\K)$~--- the class of algebras in which all finitely
generated subalgebras belong to $\K$;\\
$\pvar (\K)$~--- the least prevariety including $\K$;\\
$\qvar (\K)$~--- the least quasi-variety including $\K$, i.e.,
$\qvar (\K)=\Mod (\Th _{\rm qi} (\K))$;\\
$\ucl (\K)$~--- the universal class of algebras generated by $\K$,
i.e., $\ucl (\K)=\Mod (\Th _{\forall} (\K))$; \\$\Res (\K)$~---
the class of algebras which are separated by $\K$;
\\$\Dis (\K)$~--- the class of algebras which are discriminated by
$\K$;
\\
$\K\e$~--- the addition of the trivial algebra $\E$ to $\K$, i.e., $\K\e=\K\cup \{\E\}$;\\
$\K_\omega$~--- the class of finitely generated algebras from
$\K$.

Here we denote by $\Th_{\rm qi}(\K)$ (correspondingly,
$\Th_{\forall}(\K)$, $\Th_{\exists}(\K)$) the set of all
quasi-identities (correspondingly, universal sentences,
existential sentences) which are true in all structures from $\K$.

For an arbitrary class $\K$ of $\L$-algebras one has:
$$
\begin{array}{rl}
\ucl(\K)=\S\Pu(\K), &\; \Dis(\K) \subseteq \ucl(\K),\\
\Res(\K)=\pvar(\K)=\S\P(\K), &\; \pvar(\K)\subseteq\qvar(\K).
\end{array}
$$

According to Gorbunov~\cite{Gorbunov} and in contrast
to~\cite{DMR1}, we assume that the direct product for the empty
set of indexes coincides with the trivial $\L$-algebra $\E$. In
particularly, when we say that an algebra $\Ce$ is a finite direct
product of algebras from $\K$ (or a subdirect product of a
finitely many algebras from $\K$) then $\Ce$ may be just the
trivial algebra. However, while defining an filterproduct we
assume that the set of indexes is non-empty.

\section{Generalizations of the Unification Theorems} \label{sec:unification_theorems_compact}

Unification Theorems~A and~C are formulated in Introduction above
for an equationally Noetherian algebra $\B$. Those theorems have
been proven in~\cite{DMR1, DMR2}.

\medskip

\noindent{\bf Question:} {\em Suppose that the algebra $\B$ is not
equationally Noetherian. When Unification Theorems remain true for
$\B$?}

To answer this question we need to analyze the proofs of
Theorems~A and~C. As it was mentioned in~\cite{DMR2}, for the
reasoning of some implications in Theorems~A and~C the
equationally Noetherian property is not required, namely, one has
the following remark.

\begin{remark}\label{rem}
Let $\B$ be an algebra in a functional language $\L$ and $\Ce$ a
finitely generated $\L$-algebra. Then
\begin{itemize}
\item $\Ce$ is the coordinate algebra of an
irreducible algebraic set over $\B$ defined by a system of
equations in the language $\L$ \; IF AND ONLY IF \; $\Ce$ is
discriminated by $\B$ (Theorem~A: $7 \Longleftrightarrow 4$);
\item IF \; $\Ce$ is
discriminated by $\B$ \; THEN \; $\Ce$ is a limit algebra over
$\B$ (Theorem~A: $4 \Longrightarrow 5$);
\item $\Ce$ is the coordinate algebra of an algebraic set over $\B$
defined by a system of equations in the language $\L$ \; IF AND
ONLY IF \; $\Ce \in \pvar(\B)$ (Theorem~C: $7 \Longleftrightarrow
2$);
\item IF \; $\Ce$ is a subdirect product of a finitely many limit algebras over $\B$ \; THEN \; $\Ce \in \qvar(\B)$
(Theorem~C: $5 \Longrightarrow 1$); and so on.
\end{itemize}
The complete set of implications in Theorems~A and~C which always
remain true is represented as follows:
$$
\mbox{Theorem~A:} \quad \{4 \Leftrightarrow 7\}\quad
\Longrightarrow \quad \{1 \Leftrightarrow 2\Leftrightarrow
3\Leftrightarrow 5\Leftrightarrow 6\};
$$
$$
\mbox{Theorem~C:} \quad \{5\} \quad \Longrightarrow \quad \{1
\Leftrightarrow 6\} \quad \Longleftarrow \quad \{2 \Leftrightarrow
3 \Leftrightarrow 4 \Leftrightarrow 7\}.
$$
\end{remark}

Further, when proving $1) \Longrightarrow 4)$ in both Theorems~A
and~C, we use not equationally Noetherian property itself, but
some weaker properties. What properties exactly? These are
$\uw$-compactness and $\qw$-compactness.

\begin{definition}
We say $\L$-algebra $\B$ is {\em $\qw$-compact} if for any finite
set $X$, any system of equations $S\subseteq \At_\L (X)$, and any
equation $(t_0=s_0)\in \At_\L (X)$ such that
$$
\V_\B(S) \; \subseteq \; \V_\B(t_0=s_0)
$$
there exists a finite subsystem $S_0\subseteq S$ such that
$$
\V_\B(S)\; \subseteq \; \V_\B(S_0) \; \subseteq \; \V_\B(t_0=s_0).
$$
Here the finite subsystem $S_0$ may alter depending on equation
$(t_0=s_0)$.
\end{definition}

\begin{definition}
An $\L$-algebra $\B$ is termed {\em $\uw$-compact} if for any
finite set $X$, any system of equations $S\subseteq \At_\L (X)$,
and any equations $(t_1=s_1), \ldots, (t_m=s_m)\in \At_\L (X)$
such that
$$
\V_\B(S) \; \subseteq \; \V_\B(t_1=s_1)\:\cup \:\ldots\: \cup\:
\V_\B(t_m=s_m)
$$
there exists a finite subsystem $S_0\subseteq S$ such that
$$
\V_\B(S)\; \subseteq \; \V_\B(S_0) \; \subseteq \;
\V_\B(t_1=s_1)\:\cup \:\ldots\: \cup\: \V_\B(t_m=s_m).
$$
Here the finite subsystem $S_0$ may alter depending on equations
$(t_1=s_1),\ldots,(t_m=s_m)$.
\end{definition}

It is clear that any equationally Noetherian algebra $\B$ is
$\uw$-compact, and any $\uw$-compact algebra is $\qw$-compact.

The definitions of $\uw$-compactness and $\qw$-compactness above
are given in geometric form. We know some other approaches to
these notions that will be discussed in Section~\ref{sec:compact}.
In that section will be also represented the etymology of the
notion of $\uw$($\qw$)-compactness.

\subsection{The generalization of Unification
Theorem~A}\label{subsec:A}

The significance of $\uw$-compact algebras in universal algebraic
geometry is shown in the following theorem.

\begin{theorem}[analog of Theorem~A]\label{uw}
Let $\B$ be $\uw$-compact algebra in a functional language $\L$.
Then for a finitely generated algebra $\Ce$ of $\L$ the following
conditions are equivalent:
\begin{enumerate}
\item [1)] $\Th_{\forall} (\B) \subseteq \Th_{\forall} (\Ce)$, i.e., $\Ce \in \ucl(\B)$;
\item [2)] $\Th_{\exists} (\B) \supseteq \Th_{\exists} (\Ce)$;
\item [3)] $\Ce$ embeds into an ultrapower of $\B$;
\item [4)] $\Ce$ is discriminated by $\B$;
\item [5)] $\Ce$ is a limit algebra over $\B$;
\item [6)] $\Ce$ is an algebra defined by a complete atomic type in the theory $\Th _{\forall} (\B)$ in
$\L$;
\item [7)] $\Ce$ is the coordinate algebra of an irreducible
algebraic set over $\B$ defined by a system of equations in the
language $\L$.
\end{enumerate}
Moreover, if for an $\L$-algebra $\B$ and for all finitely
generated $\L$-algebras $\Ce$ the conditions above are equivalent
then $\B$ is $\uw$-compact.
\end{theorem}

\begin{proof}
It follows from Remark~\ref{rem} that conditions 1)--7) are
equivalent if and only if one has equivalence
$1)\Longleftrightarrow 4)$. The latter means that a finitely
generated algebra $\Ce$ is discriminated by $\B$ if and only if
$\Ce \in \ucl(\B)$, i.e., $\ucl (\B)_\omega=\Dis (\B)_\omega$. By
Theorem~\ref{theoremU} below, one has the equality $\ucl
(\B)_\omega=\Dis (\B)_\omega$ if and only if an algebra $\B$ is
$\uw$-compact.
\end{proof}

\subsection{The generalization of Unification
Theorem~C}\label{subsec:C}

To prove an analog of Theorem~C for $\qw$-compact algebras we need
the following results.

\begin{lemma}[\cite{DMR1}]\label{Corollary 5.7}
Let $\Ce$ be a limit algebra over an $\L$-algebra $\B$. Then there
exists an ultrapower $\B^\ast$ of $\B$ such that $\Ce$ embeds into
$\B^\ast$.
\end{lemma}

\begin{lemma}[\cite{DMR2}]\label{Lemma 3.9}
A finitely generated $\L$-algebra $\Ce$ is the coordinate algebra
of an algebraic set over $\L$-algebra $\B$ if and only if $\Ce$ is
a subdirect product of the coordinate algebras of irreducible
algebraic sets over $\B$.
\end{lemma}

\begin{theorem}[analog of Theorem~C]\label{qw}
Let $\B$ be $\qw$-compact algebra in a functional language $\L$.
Then for a finitely generated algebra  $\Ce$ of $\L$ the following
conditions are equivalent:
\begin{enumerate}
\item [1)] $\Ce \in \qvar(\B)$, i.e., $\Th_{\rm qi} (\B) \subseteq \Th_{\rm qi} (\Ce)$;
\item [2)] $\Ce \in \pvar(\B)$;
\item [3)] $\Ce$ embeds into a direct power of $\B$;
\item [4)] $\Ce$ is separated by $\B$;
\item [5')] $\Ce$ is a subdirect product of limit algebras over $\B$;
\item [6)] $\Ce$ is an algebra defined by a complete atomic type in the theory $\Th _{\rm qi} (\B)$ in $\L$;
\item [7)] $\Ce$ is the coordinate algebra of an algebraic set over $\B$ defined by a system of
equations in the language $\L$.
\end{enumerate}
Moreover, if for an $\L$-algebra $\B$ and for all finitely
generated $\L$-algebras $\Ce$ the conditions above are equivalent
then $\B$ is $\qw$-compact.
\end{theorem}

\begin{proof}
By Remark~\ref{rem}, it is sufficient to prove implications $1)
\Longrightarrow 2)$, $5') \Longrightarrow 1)$, and $7)
\Longrightarrow 5')$ for $\qw$-compact algebra $\B$. By
Theorem~\ref{theoremQ} below, we have the identity
$\qvar(\B)_\omega=\pvar(\B)_\omega$ that gives proof of  $1)
\Longrightarrow 2)$. For implication $5') \Longrightarrow 1)$ we
refer to Lemma~\ref{Corollary 5.7} and the fact that every
quasi-variety is closed under ultraproducts, direct products and
subalgebras.

For proving $7) \Longrightarrow 5')$ suppose that $\Ce$ is the
coordinate algebra of an algebraic set over $\B$. By
Lemma~\ref{Lemma 3.9}, $\Ce$ is a subdirect product of coordinate
algebras of irreducible algebraic sets over $\B$. By
Remark~\ref{rem} (Theorem~A: $7 \Longrightarrow 5$), coordinate
algebras of irreducible algebraic sets over $\B$ are limit
algebras over $\B$.

Suppose now that for some $\L$-algebra $\B$ we have equivalence
$1) \Longleftrightarrow 2)$ for all finitely generated
$\L$-algebras $\Ce$. It means that
$\qvar(\B)_\omega=\pvar(\B)_\omega$ and, by Theorem~\ref{theoremQ}
below, the algebra $\B$ is $\qw$-compact.
\end{proof}

\begin{remark}\label{problem}
Unfortunately, we are not in a position to formulate Theorem~C for
$\qw$-compact algebras in all its fullness, because item 5)
essentially needs equationally Noetherian property. We have to
weak 5), namely we should erase words ``finitely many''.
\end{remark}

To establish Remark~\ref{problem} we formulate the following
problem.

\medskip

\noindent {\bf Embedding Problem.} {\it Let $\B$ be $\qw$-compact
algebra in a functional language $\L$. The question: whether or
not every coordinate algebra over $\B$ subdirectly embeds into a
finite direct product of algebras from $\ucl (\B)$? If the answer
is ``not'', then we ask whether or not the same holds for at least
$\uw$-compact algebras.}

A.\,N.\,Shevlyakov in~\cite{Shevl1} gives the negative answer to
the Embedding Problem both for $\qw$-compact and $\uw$-compact
algebras.

Let us put an addition to Remark~\ref{rem}.

\begin{remark}
The following implications and equivalencies from Theorem~\ref{qw}
hold for an arbitrary algebra $\B$:
$$
\xymatrix{ \{1 \Leftrightarrow 6\} && \{ 2 \Leftrightarrow 3 \Leftrightarrow 4 \Leftrightarrow 7\} \ar@{=>}[ll] \ar@/^/@{=>}[dl]\\
& \{5'\} \ar@/^/@{=>}[ul]}
$$
\end{remark}

Theorem~\ref{qw} gives a classification of coordinate algebras in
terms of quasivarieties. Thereby, any characterizations of
quasivariety $\qvar (\K)$ of a class $\K$ of $\L$-algebras are
helpful in universal algebraic geometry. In~\cite{Gorbunov,
Malcev} one can find the identities:
\begin{gather*}
\qvar (\K)\:=\:\S \Pf (\K)\e\:=\:\S \P \Pu (\K) \:=\: \S \Pu \P
(\K)\:=\:\S \Pu \Pw (\K)\:=\:\\\:=\:\S \Ls \P (\K)\:=\:\Ls \S \P
(\K)\:=\:\Ls \Ps (\K)\:=\:\Ld \S \P (\K).
\end{gather*}

\subsection{Weak generalization of Unification
Theorem~A}\label{subsec:weak}

Let $\B$ be an algebra in a functional language $\L$. Let us
consider the class $\ucl (\B)_\omega$.

By Remark~\ref{rem}, for any irreducible algebraic set $Y$ over
$\B$ the coordinate algebra $\Gamma (Y)$ belongs to $\ucl
(\B)_\omega$. If $\B$ is $\uw$-compact algebra then, by
Theorem~\ref{uw}, every algebra $\Ce$ from $\ucl (\B)_\omega$ is
the coordinate algebra of some irreducible algebraic set $Y$ over
$\B$.

Let us apply a weak mode to $\uw$-compactness and require that
every coordinate algebra $\Ce$ from $\ucl (\B)_\omega$ is
irreducible. Suppose that some algebras from $\ucl (\B)_\omega$
are not coordinate algebras for algebraic sets over $\B$ at all,
however, if $\Gamma (Y) \in \ucl(\B)$ then $Y$ is irreducible. Let
us introduce a specific name for algebra $\B$ with this type of
property.

\begin{definition}
We name an $\L$-algebra $\B$ {\em weakly $\uw$-compact} if each
non-empty algebraic set $Y$ over $\B$ which coordinate algebra
$\Gamma (Y)$ belongs to $\ucl (\B)$ is irreducible.
\end{definition}

By Theorem~\ref{uw}, every $\uw$-compact algebra is weakly
$\uw$-compact. We will discuss weakly $\uw$-compact algebras,
their properties and equivalent approaches to them in
Subsection~\ref{subsec:weakcompact}.

For weakly $\uw$-compact algebras we have just the following weak
analog of Theorem~A. It allows to describe irreducible coordinate
algebras inside the class of all coordinate algebras.

\begin{theorem}[weak analog of Theorem~A]\label{wuw}
Let $\B$ be a weakly $\uw$-compact algebra in a functional
language $\L$ and $Y$ a non-empty algebraic set over $\B$. Then
the following conditions are equivalent:
\begin{enumerate}
\item [1)] $\Th_{\forall} (\B) \subseteq \Th_{\forall} (\Gamma(Y))$, i.e., $\Gamma(Y) \in \ucl(\B)$;
\item [2)] $\Th_{\exists} (\B) \supseteq \Th_{\exists} (\Gamma(Y))$;
\item [3)] $\Gamma(Y)$ embeds into an ultrapower of $\B$;
\item [4)] $\Gamma(Y)$ is discriminated by $\B$;
\item [5)] $\Gamma(Y)$ is a limit algebra over $\B$;
\item [6)] $\Gamma(Y)$ is an algebra defined by a complete atomic type in the theory $\Th _{\forall} (\B)$ in
$\L$;
\item [7)] $Y$ is irreducible.
\end{enumerate}
Moreover, if for an $\L$-algebra $\B$ and for every non-empty
algebraic set $Y$ the conditions above are equivalent then $\B$ is
weakly $\uw$-compact.
\end{theorem}

\begin{proof}
It follows from Remark~\ref{rem} that conditions~1)--7) are
equivalent if and only if one has implication $1) \Longrightarrow
7)$. By definition, implication $1) \Longrightarrow 7)$ take place
if and only if $\B$ is weakly $\uw$-compact.
\end{proof}

\section{$\qw$-compact and $\uw$-compact algebras}
\label{sec:compact}

In Section~\ref{sec:unification_theorems_compact} we gave the
definitions of $\qw$- and $\uw$-compact algebras in geometric
language. In Subsection~\ref{subsec:compact} we gather the
numerous another approaches to these notions into two theorems. We
will prove these theorems in Subsection~\ref{subsec:proof}.

In Subsection~\ref{subsec:local} we introduce ``local
$\qw$($\uw$)-compact property'' and show its use in universal
algebraic geometry. Subsection~\ref{subsec:Ecompact} contains some
accessory materials.

\subsection{Criteria of $\qw$- and $\uw$-compactness} \label{subsec:compact}

At first we formulate the theorems and then give the necessary
explanations.

\begin{theorem}\label{theoremQ}
For an algebra $\B$ in a functional language $\L$ the following
conditions are equivalent:
\begin{itemize}
\item[1)] $\B$ is $\qw$-compact;
\item[2)] for any finite set $X$, any system of equations $S\subseteq \At_\L
(X)$,
and any consequence $c=(t_0=s_0)\in \Rad_\B (S)$ there exists a
finite subsystem $S_c\subseteq S$ such that $c\in \Rad_\B(S_c)$;
\item[3)] for any finite set $X$, any subset $S\subseteq \At_\L (X)$, and any atomic
formula $(t_0=s_0)\in \At_\L (X)$ if an (infinite) formula
$$
\forall \;  x_1 \ldots \forall \;  x_n \left( \bigwedge
\limits_{(t=s) \in S} t (\bar{x})=s (\bar{x}) \; \longrightarrow
\; t_0 (\bar{x})=s_0(\bar{x}) \right)
$$
holds in $\B$ then for some finite subsystem $S_c\subseteq S$ the
quasi-identity
$$
\forall \;  x_1 \ldots \forall \;  x_n \left( \bigwedge
\limits_{(t=s) \in S_c} t (\bar{x})=s (\bar{x}) \; \longrightarrow
\; t_0 (\bar{x})=s_0(\bar{x}) \right)
$$
also holds in $\B$;
\item[4)] for any finite set $X$, any subset $S\subseteq \At_\L (X)$, and any atomic
formula $(t_0=s_0)\in \At_\L (X)$ if the set of formulas
$$
T\; = \; S\cup \{\neg (t_0=s_0)\}
$$
is finitely satisfiable in $\B$ then it is satisfiable in $\B$;
\item[5)] every finitely generated algebra from $\qvar(\B)$ is the
coordinate algebra of an algebraic set over $\B$;
\item[6)] $\qvar(\B)_\omega=\pvar(\B)_\omega$;
\item[7)] $\qvar(\B)=\Lo\Res (\B)$;
\item[8)] $\Ls \S \P (\B)=\Lo \S \P (\B)$;
\item[9)] $\Ld \S \P (\B)=\Lo \S \P (\B)$;
\item[10)] for any
finite set $X$ and any system of equations $S\subseteq \At_\L (X)$
one has:
$$
\Rad_\B(S)=\bigcup_{S_0 \subseteq S} \: \Rad _\B (S_0),
$$
where $S_0$ runs all finite subsystems of $S$;
\item[11)] for any finite set $X$ and any directed system $\{S_i, i\in I\}$
of radical ideals over $\B$ from $\At_\L(X)$ the union
$S=\bigcup_{i\in I}{S_i}$ is a radical ideal over $\B$;
\item[12)] for any finite set $X$ and any epimorphic direct system
$\Lambda= (I,\Ce_i,h_{ij})$ of coordinate algebras over $\B$ with
generating set $X$, and $h_{ij}(x)=x$, $x\in X$, the epimorphic
direct limit $\underrightarrow{\lim} \: \Ce_i$ is a coordinate
algebra over $\B$.
\end{itemize}
\end{theorem}

\begin{theorem}\label{theoremU}
For an algebra $\B$ in a functional language $\L$ the following
conditions are equivalent:
\begin{itemize}
\item[1)] $\B$ is $\uw$-compact;
\item[2)] for any
finite set $X$, any system of equations $S\subseteq \At_\L (X)$,
and any consequence $c$ of $S$ over $\B$ of the form
$c=(t_1=s_1)\vee\ldots \vee (t_m=s_m)$, $t_i,s_i\in \Tr_\L(X)$,
there exists a finite subsystem $S_c\subseteq S$ such that $c$ is
a consequence of $S_c$ over $\B$;
\item[3)] for any finite set $X$, any subset $S\subseteq \At_\L
(X)$,
and any atomic formulas $(t_1=s_1), \ldots, (t_m=s_m)\in \At_\L
(X)$ if an (infinite) formula
$$
\forall \;  x_1 \ldots \forall \;  x_n \left( \bigwedge
\limits_{(t=s) \in S} t (\bar{x})=s (\bar{x}) \; \longrightarrow
\; \bigvee\limits_{i=1}^m t_i (\bar{x})=s_i(\bar{x}) \right)
$$
holds in $\B$ then for some finite subsystem $S_c\subseteq S$ the
universal sentence
$$
\forall \;  x_1 \ldots \forall \;  x_n \left( \bigwedge
\limits_{(t=s) \in S_c} t (\bar{x})=s (\bar{x}) \; \longrightarrow
\; \bigvee\limits_{i=1}^m t_i (\bar{x})=s_i(\bar{x}) \right)
$$
also holds in $\B$;
\item[4)] for any finite set $X$, any subset $S\subseteq \At_\L (X)$, and
any atomic formulas $(t_1=s_1), \ldots, (t_m=s_m)\in \At_\L (X)$
if the set of formulas
$$
T\; =\; S\cup \{\neg (t_1=s_1), \ldots , \neg (t_m=s_m)\}
$$
is finitely satisfiable in $\B$ then it is satisfiable in $\B$;
\item[5)] every finitely generated algebra from $\ucl(\B)$ is the
coordinate algebra of an irreducible algebraic set over $\B$;
\item[6)] $\ucl(\B)_\omega=\Dis(\B)_\omega$;
\item[7)] $\ucl(\B)=\Lo \Dis (\B)$.
\end{itemize}
\end{theorem}

Item~2) in Theorem~\ref{theoremQ} (correspondingly, in
Theorem~\ref{theoremU}) gives the definition of $\qw$-compact
(correspondingly, $\uw$-compact) algebra in terms of radicals;
item~3)~--- in terms of infinite formulas; item~5)~--- in terms of
coordinate algebras.

Item~4) shows that the definition of $\qw$($\uw$)-compactness is a
compact property relating to special types of sets of formulas
$T$, as it is discussed in Introduction. The background of this
notion is detailed in~\cite{MR2} for groups. Here we will tell
just a few words about it.

The answer for the following question has been attained by
V.\,A.\,Gorbunov~\cite{Gorbunov}.

\medskip

\noindent {\bf Malcev Problem.} {\it When the prevariety
$\pvar(\K)$ generated by class $\K$ is a quasivariety?}

V.\,A.\,Gorbunov has introduced the notion of quasi-compact (${\rm
q}$-compact) class $\K$ and proved that $\pvar (\K)=\qvar(\K)$ if
and only if $\K$ is ${\rm q}$-compact. Let us compare that result
with item~6) in Theorem~\ref{theoremQ}.

The definition of ${\rm q}$-compact algebra $\B$ is much the same
as the definition of $\qw$-compact algebra given in item~4) of
Theorem~\ref{theoremQ}. We just bound the set of variables $X$ for
defining $\qw$-compact algebras: $X$ must be finite. For ${\rm
q}$-compact algebras $X$ runs sets of all possible cardinalities.

While items~1)--7) in Theorems~\ref{theoremQ} and~\ref{theoremU}
are symmetric, items~10)--12) in Theorem~\ref{theoremQ} are
specific for $\qw$-compact algebras; 8) and~9) in
Theorem~\ref{theoremQ} are just corollaries of~7).

Items~10) and~11) in Theorems~\ref{theoremQ} are close. The family
$\{\Rad (S_0)\}$, where $S_0$ runs all finite subsystems of a
system $S$, gives an example of a directed system. Let us remind
concerned definitions.

A partial ordering  $(I, \leqslant)$ is {\em directed} if any two
elements from $I$ have an upper bound. A family $\{\theta_i, i\in
I\}$ of congruencies on an $\L$-algebra $\M$ with $i\leqslant
j\Leftrightarrow\theta_i\subseteq\theta_j$ is called {\em directed
system of congruencies}.

A system $S \subseteq \At_\L(X)$ is {\em radical ideal} over $\B$
if $S=\Rad_\B (S)$.

\begin{definition}
We say that a family $\{S_i, i\in I\}$ of radical ideals from
$\At_\L(X)$ is a {\em directed system} if the family
$\{\theta_{S_i}, i\in I\}$ is a directed system of congruencies on
$\T_\L(X)$.
\end{definition}

Let us prove just a little part of Theorems~\ref{theoremQ}
and~\ref{theoremU}.

\begin{lemma}\label{lemmacompact}
Let $\B$ be an $\L$-algebra, $X$ a finite set, $|X|=n$,
$S\subseteq\At_\L (X)$ a system of equations, and $(t_1=s_1),
\ldots, (t_m=s_m)\in \At_\L (X)$ atomic formulas. Then the
following conditions are equivalent:
\begin{itemize}
\item[1)] $\V_\B(S) \; \subseteq \; \V_\B(t_1=s_1)\:\cup \:\ldots\:
\cup\: \V_\B(t_m=s_m)$;
\item[2)]  $(t_1=s_1)\vee\ldots \vee
(t_m=s_m)$ is a consequence of $S$ over $\B$;
\item[3)] the (infinite) formula
$$
\forall \;  x_1 \ldots \forall \;  x_n \left( \bigwedge
\limits_{(t=s) \in S} t (\bar{x})=s (\bar{x}) \; \longrightarrow
\; \bigvee\limits_{i=1}^m t_i (\bar{x})=s_i(\bar{x}) \right)
$$
holds in $\B$;
\item[4)] the set of formulas
$$
T\; =\; S\cup \{\neg (t_1=s_1), \ldots , \neg (t_m=s_m)\}
$$
is not satisfiable in $\B$;
\item[5)] there is no homomorphism $h\colon \langle \{c_1,\ldots,c_n \} \,|\, S\rangle \to
\B$ such that
$$
h(t_i(c_1,\ldots,c_n))\ne h(s_i(c_1,\ldots,c_n)) \quad \mbox{for
all} \quad i \in \{1,\ldots, m\}.
$$
\end{itemize}
\end{lemma}

\begin{proof}
Straightforward.
\end{proof}

\begin{corollary}
One has equivalencies $1) \Longleftrightarrow 2)$, $1)
\Longleftrightarrow 3)$, $3) \Longleftrightarrow 4)$ in both
Theorems~\ref{theoremQ} and~\ref{theoremU}.
\end{corollary}

\begin{proof}
Equivalencies $1) \Longleftrightarrow 2)$, $1) \Longleftrightarrow
3)$ are easy. Note that the statement in item~3) has a form ``$A$
implies $B$''. The equivalent statement is ``$\neg B$ implies
$\neg A$'' which gives~4). So we have $3) \Longleftrightarrow 4)$.
\end{proof}

From now on, we will use not only geometric definition of
$\qw$-compact (correspondingly, $\uw$-compact) algebra, but also
the definitions that items~2),~3),~4) in Theorem~\ref{theoremQ}
(correspondingly, in Theorem~\ref{theoremU}) give us.

\subsection{$\E$-compact algebras}\label{subsec:Ecompact}

This subsection is a special excursus. We consider here the
following problem.

\medskip

\noindent {\bf Problem.} {\it When the conditions ``$\E\in
\ucl(\B)$'' and ``$\B$ has a trivial subalgebra'' are equivalent?}

It is important to note that for a large class of algebras the
conditions ``$\E\in \ucl(\B)$'' and ``$\B$ has a trivial
subalgebra'' are equivalent, but not for all algebras.

\begin{definition}
We say an $\L$-algebra $\B$ is {\em $\E$-compact} if finite
satisfiability in $\B$ of the set of all atomic formulas
$\At_\L(\{x\})$ in one variable $x$ implies its satisfiability in
$\B$.
\end{definition}

\begin{lemma}\label{lemma2}
An $\L$-algebra $\B$ is $\E$-compact if and only if the conditions
``$\E\in \ucl(\B)$'' and ``$\B$ has a trivial subalgebra'' are
equivalent.
\end{lemma}

\begin{proof}
It is sufficient to show that $\At_\L(\{x\})$ is satisfiable in
$\B$ if and only if $\B$ has a trivial subalgebra, and
$\At_\L(\{x\})$ is finitely satisfiable if and only if $\E \in
\ucl (\B)$.

Suppose that $\At_\L(\{x\})$ is satisfiable in $\B$. Then there
exists an element $b\in B$ with $\B \models (t(b)=s(b))$ for all
$t,s\in\Tr_\L(\{x\})$. Therefore, subalgebra of $\B$ generated by
the element $b$ is trivial. Conversely, if $\B$ has a trivial
subalgebra $\E=\{e\}$ then the set of all atomic formulas
$\At_\L(\{x\})$ is realized in $\B$ on the element $e$.

Assume now that $\At_\L(\{x\})$ is not finitely satisfiable in
$\B$. Then there exists a finite set $S_0$ of atomic formulas such
that the universal sentence
\begin{equation}\label{eq:e}
\forall \; x\quad \left(\bigvee_{(t=s)\in S_0} {\neg(\:t(x)=
s(x)\:)} \right)
\end{equation}
holds in $\B$. However~\eqref{eq:e} is false in $\E$, so $\E
\not\in \ucl (\B)$. Conversely, if the set of all atomic formulas
$\At_\L(\{x\})$ is finitely satisfiable in $\B$ then by
Compactness Theorem it is realized in some ultrapower $\B^\ast$ of
$\B$. Hence, $\E \in \ucl (\B)$.
\end{proof}

\begin{corollary}
The condition ``algebra $\B$ is $\E$-compact'' means that $\B$ has
a trivial subalgebra or $\E\not\in\ucl(\B)$.
\end{corollary}

Let us note that in ``good'' signatures all algebras are
$\E$-compact.

\begin{lemma}\label{lemmaL}
Suppose a functional language $\L$ has at least one constant
symbol. Then every algebra in $\L$ is $\E$-compact.
\end{lemma}

\begin{proof}
Let $\B$ be an $\L$-algebra. We need to show that condition $\E\in
\ucl(\B)$ implies that $\B$ has a trivial subalgebra. Consider the
set of formulas
$$
T=\{c = c'\} \cup \{F(c,\ldots,c)=c\},
$$
where $c,c'$ run all constant symbols from $\L$ and $F$ runs all
functional symbols from $\L$. If $\E\in \ucl(\B)$, then $\B\models
T$. Therefore, there exists an element $b\in B$ such that $c^\B=b$
for all constant symbol $c$ from $\L$, and $F(b,\ldots,b)=b$ for
all functional symbol $F$ from $\L$. Thereby, the element $b$
generates the trivial subalgebra in $\B$.
\end{proof}

\begin{lemma}
Suppose $\L$ is a finite functional language. Then every algebra
in $\L$ is $\E$-compact.
\end{lemma}

\begin{proof}
After Lemma~\ref{lemmaL} we may assume that $\L$ has no constant
symbols. Let $\B$ an $\L$-algebra. If $\E\in \ucl(\B)$ then the
existential sentence
$$
\exists \; x\quad \left(\bigwedge_{F\in \: \L} {F(x,\ldots,x)=x}
\right)
$$
holds in $\B$. Thereby, $\B$ has a trivial subalgebra.
\end{proof}

If $\L$ is an infinite functional language with no constant
symbols, then it is easy to construct an $\L$-algebra $\B$ that is
not $\E$-compact (see Example~\ref{example} below).

It follows from the definition that all equationally Noetherian
algebras are $\E$-compact. Now we state that all $\qw$- and
$\uw$-compact algebras are $\E$-compact. We need the following
facts and definitions.

According to V.\,A.\,Gorbunov~\cite{Gorbunov}, an $\L$-algebra
$\B$ is {\em weakly atomic compact}, if for any set $X$ and any
subset $S \subseteq \At_\L(X)$ finite satisfiability of $S$ in
$\B$ implies realizability of $S$ in $\B$. We say that an
$\L$-algebra $\B$ is {\em weakly atomic $\omega$-compact}, if for
any finite set $X$ and any subset $S \subseteq \At_\L(X)$ finite
satisfiability of $S$ in $\B$ implies realizability of $S$ in
$\B$. It is obvious that weak atomic $\omega$-compactness implies
$\E$-compactness.

The following result has been proven by M.\,Kotov~\cite{Kotov}.

\noindent{\bf Lemma} (\cite{Kotov}){\bf.} {\it Every $\qw$-compact
algebra in a functional language $\L$ is weakly atomic
$\omega$-compact.}

\begin{corollary}\label{corkot}
Let $\B$ be $\qw$-compact $\L$-algebra (in particularly, $\B$ may
be $\uw$-compact). Then the universal closure $\ucl (\B)$ contains
the trivial algebra $\E$ if and only if $\B$ has a trivial
subalgebra.
\end{corollary}

Let us note that M.\,Kotov has proven more general result in his
work. We formulate it on geometric language.

\noindent{\bf Lemma} (\cite{Kotov}){\bf.} {\it Let $\B$ be an
$\L$-algebra and $S$ a system of equations in $\L$. If $\B$ is
$\qw$-compact and $\V_\B(S)$ is a singleton set or the empty set,
then there exists a finite subsystem $S_0\subseteq S$ which is
equivalent to $S$ over $\B$. If $\B$ is $\uw$-compact and
$\V_\B(S)$ is a finite set or the empty set, then there exists a
finite subsystem $S_0\subseteq S$ which is equivalent to $S$ over
$\B$.}

\subsection{Local compact properties} \label{subsec:local}

Let $X$ be a finite set. Fix a subset $S\subseteq\At_\L (X)$. We
will give the definitions of local compact properties with respect
to fixed $S$.

\begin{definition}
An $\L$-algebra $\B$ is called {\em $\qs$-compact} if for each
atomic formula $(t_0=s_0)\in \At_\L (X)$ if the set of formulas
$$
T\; = \; S\cup \{\neg (t_0=s_0)\}
$$
is finitely satisfiable in $\B$ then it is satisfiable in $\B$.
\end{definition}

\begin{definition}
An $\L$-algebra $\B$ is called {\em $\us$-compact} if for any
atomic formulas $(t_1=s_1), \ldots, (t_m=s_m)\in \At_\L (X)$ if
the set of formulas
\begin{equation}\label{eq:T}
T\; =\; S\cup \{\neg (t_1=s_1), \ldots , \neg (t_m=s_m)\}
\end{equation}
is finitely satisfiable in $\B$ then it is satisfiable in $\B$.
\end{definition}

It is clear that algebra $\B$ is $\qw$($\uw$)-compact if and only
if it is $\qs$($\us$)-compact for every finite set $X$ and every
$S\subseteq\At_\L (X)$.

The main results on local compact properties are the following.

\begin{proposition}\label{localQ}
Let $\B$ be an algebra in a functional language $\L$, $X$ a finite
set, $S\subseteq\At_\L (X)$, and $\Ce=\langle X | \, S\rangle$.
Then the following conditions are equivalent:
\begin{enumerate}
\item [1)] $\Ce$ is the coordinate algebra of an
algebraic set over $\B$ defined by a system of equations in the
language $\L$;
\item [2)] $\Ce$ is separated by $\B$;
\item [3)] $\Ce\in \qvar(\B)$ and $\B$ is $\qs$-compact.
\end{enumerate}
\end{proposition}

\begin{proposition}\label{localU}
Let $\B$ be an algebra in a functional language $\L$, $X$ a finite
set, $S\subseteq\At_\L (X)$, such that $[S]\ne \At_\L (X)$, and
$\Ce=\langle X | \, S\rangle$. Then the following conditions are
equivalent:
\begin{enumerate}
\item [1)] $\Ce$ is the coordinate algebra of an irreducible
algebraic set over $\B$ defined by a system of equations in the
language $\L$;
\item [2)] $\Ce$ is discriminated by $\B$;
\item [3)] $\Ce\in\ucl(\B)$ and $\B$ is $\us$-compact.
\end{enumerate}
\end{proposition}

Before giving a proof of these propositions we need some remarks.
Firstly, equivalence $1) \Longleftrightarrow 2)$ in both
Propositions~\ref{localQ} and~\ref{localU} have been proven
in~\cite{DMR2}. Secondly, let us answer the question: when the set
of formulas~\eqref{eq:T} is not finitely satisfiable in $\B$? It
happens if and only if there exists a finite subset $S_0\subseteq
S$ such that the universal sentence
\begin{equation}\label{eq:unfor}
\forall \;  y_1 \ldots \forall \;  y_n \left( \bigwedge
\limits_{(t=s) \in S_0} t (\bar{y})=s (\bar{y}) \; \longrightarrow
\; \bigvee \limits_{i=1}^{m} t_i (\bar{x})=s_i
(\bar{y})\right),\quad\mbox{where} \; |X|=n,
\end{equation}
holds in $\B$. For example, if $(t_i=s_i)\in [S]$ for some
$i\in\{1,\ldots,m\}$, then there exists a finite subset
$S_0\subseteq S$ such that $S_0\vdash (t_i=s_i)$, in particularly,
universal formula~\eqref{eq:unfor} holds in $\B$.

Thirdly, note that in Propositions~\ref{localU} we claim $[S]\ne
\At_\L (X)$, but in Propositions~\ref{localQ} such restriction is
omitted. If $[S]=\At_\L (X)$ then $\Ce=\langle X | \, S\rangle$ is
the trivial algebra $\E$. Moreover, in this case every algebra
$\B$ is $\qs$- and $\us$-compact. Since the trivial algebra $\E$
is the coordinate algebra of an algebraic set over $\B$ anyway and
$\E$ belongs to each quasi-variety~\cite{DMR2}, we have no
difficulties with $\E$ in Propositions~\ref{localQ}.

\begin{remark}
One can omit restriction $[S]\ne\At_\L(X)$ in
Proposition~\ref{localU} if and only if $\B$ is $\E$-compact
algebra. Indeed, the trivial algebra $\E$ is the coordinate
algebra of an irreducible algebraic set over $\B$ if and only if
$\B$ has a trivial subalgebra~\cite[Lemma~3.22]{DMR2}. By
Lemma~\ref{lemma2}, the conditions ``$\E\in \ucl(\B)$'' and ``$\B$
has a trivial subalgebra'' are equivalent if and only if $\B$ is
$\E$-compact.
\end{remark}

Now we are going to prove Propositions~\ref{localQ}
and~\ref{localU}. Arguments for them are the similar, so we will
prove only Propositions~\ref{localU}.

\begin{proof}[Proof of Propositions~\ref{localU}]
Let $\Ce\simeq \T_{\L}(X)/\theta_S$, $X=\{c_1,\ldots,c_n\}$, and
$[S]\ne\At_\L (X)$. By definition $\Ce$ is discriminated by $\B$
if for any finite set of atomic formulas $(t_1=s_1), \ldots,
(t_m=s_m)\in\At_\L(X)\setminus [S]$ there exists a homomorphism
$h\colon \Ce \to \B$, such that $h(t_i(c_1,\ldots,c_n))\ne
h(s_i(c_1,\ldots,c_n))$ for all $i\in\{1,\ldots,m\}$. The
existence of such homomorphism $h\colon \Ce \to \B$ means that the
set $T$ in~\eqref{eq:T} is realized in $\B$. Note that if we take
$(t_i=s_i)\in [S]$ for some $i\in\{1,\ldots,m\}$, then $T$ is not
finitely satisfiable in $\B$. Anyway, we shown that if $\Ce$ is
discriminated by $\B$ then $\B$ is $\us$-compact. The occurrence
$\Ce\in\ucl(\B)$ follows from the inclusion $\Dis (\B) \subseteq
\ucl(\B)$.

Suppose now that $\Ce=\langle X | \, S\rangle$ is not
discriminated by $\B$ and show that $\Ce\not\in\ucl(\B)$ or $\B$
is not $\us$-compact. In this case for some atomic formulas
$(t_1=s_1), \ldots, (t_m=s_m)\in\At_\L(X)\setminus [S]$ the set
$T$ from~\eqref{eq:T} is not realized in $\B$. If at the same time
$T$ is finitely satisfiable in $\B$ then $\B$ is not
$\us$-compact. Assume that $T$ is not finitely satisfiable in
$\B$. Therefore, there exists a finite subset $S_0\subseteq S$
such that the universal formula~\eqref{eq:unfor} holds in $\B$. On
the other hand, the formula
$$
\bigwedge \limits_{(t=s) \in S_0} t (\bar{y})=s (\bar{y}) \;
\longrightarrow \; \bigvee \limits_{i=1}^{m} t_i (\bar{x})=s_i
(\bar{y})
$$
is false in $\Ce$ under the interpretation $y_i\mapsto c_i$, $i=1,
\ldots,n$, hence  $\Ce \not\in \ucl (\B)$.
\end{proof}

\subsection{Proof of the criteria}\label{subsec:proof}

In this subsection we prove Theorems~\ref{theoremQ}
and~\ref{theoremU} that have been formulated in
Subsection~\ref{subsec:compact}. Remain that equivalencies $1)
\Longleftrightarrow 2)$, $1) \Longleftrightarrow 3)$, $3)
\Longleftrightarrow 4)$ in both theorems have been proven in
Subsection~\ref{subsec:compact}.

At first we prove the following easy lemma that will be useful
below.

\begin{lemma}\label{lemmaT0}
Let $\B, \Ce$ be $\L$-algebras, $\Ce \in \ucl(\B)$, and $T$ a set
of quantifier-free formulas in $\L$. If $T$ is finitely
satisfiable in $\Ce$ then it is finitely satisfiable in $\B$.
\end{lemma}

\begin{proof}
Suppose $T$ is finitely satisfiable in $\Ce$. Then for every
finite subset $\{\phi_1,\ldots,\phi_m\}\subseteq T$ the
existential sentence
\begin{equation}\label{eq:phi}
\exists \;  x_1 \ldots \exists \;  x_n \quad \left(
\phi_1(x_1,\ldots,x_n) \:\wedge\: \ldots \:\wedge\:
\phi_m(x_1,\ldots,x_n) \right)
\end{equation}
holds in $\Ce$. Since $\Ce \in \ucl(\B)$ then~\eqref{eq:phi} holds
in $\B$ too. Thereby, $T$ is finitely satisfiable in $\B$.
\end{proof}

We start with Theorem~\ref{theoremU}. Consider item~6). It states
that $\ucl(\B)_\omega=\Dis(\B)_\omega$. As inclusion $\ucl
(\B)_\omega \supseteq \Dis (\B)_\omega$ holds for an arbitrary
algebra $\B$, then item~6) is equivalent to inclusion $\ucl
(\B)_\omega \subseteq \Dis (\B)_\omega$. On the other hand,
$\Dis(\B)_\omega$ is the class of all irreducible coordinate
algebras over $\B$~\cite[Corollary 3.39]{DMR2}. Hence, we have
equivalence $5) \Longleftrightarrow 6)$.

Now let us show equivalence $4) \Longleftrightarrow 6)$. Suppose
$\B$ is $\uw$-compact and $\M$ is a finitely generated algebra
from $\ucl (\B)$. If $\M$ is a trivial algebra then, by
Corollary~\ref{corkot}, $\B$ has a trivial subalgebra, therefore,
$\M$ is discriminated by $\B$.

For non-trivial algebra $\M$ let us find a presentation $\langle X
\mid S\rangle$, where $X$ is a finite set and $S \subseteq
\At_{\L}(X)$, $[S] \ne \At_\L (X)$. As $\B$ is $\us$-compact we
have $\M \in \Dis(\B)$, by Proposition~\ref{localU}. Thus we
proved inclusion $\ucl (\B)_\omega \subseteq \Dis (\B)_\omega$ and
implication $4) \Longrightarrow 6)$.

We prove the converse implication $6) \Longrightarrow 4)$ by
contradiction. Suppose that $\B$ is not $\uw$-compact. Then there
exists a finite set $X$, a subset $S \subseteq \At_{\L}(X)$, and
atomic formulas $(t_1=s_1),\ldots,(t_m=s_m) \in \At_\L (X)$, such
that the set of formulas
$$
T\; =\; S\cup \{\neg (t_1=s_1), \ldots , \neg (t_m=s_m)\},
$$
is not realized in $\B$, but every its finite subset is realized
in $\B$.

By Compactness Theorem $T$ is realized in some ultrapower $\B^I/D$
of $\B$. Let $c_1, \ldots, c_n$ be elements from $\B^I/D$, such
that $\B^I/D \models T(c_1,\ldots,c_n)$, and $\Ce$ subalgebra of
$\B^I/D$ generated by the set $\{c_1,\ldots,c_n\}$. Clearly,
$\Ce$ is finitely generated algebra from $\ucl (\B)$. Show that
$\Ce$ is not discriminated by $\B$.

Let $\langle \{c_1,\ldots,c_n\} \mid R \rangle$ be a presentation
of $\Ce$, i.e., $\Ce \simeq \T_\L(X)/\theta_R$, $R\subseteq \At_\L
(X)$. Since $\Ce \models T(c_1,\ldots,c_n)$, one has $S\subseteq
R$ and $(t_i=s_i)\not\in [R]$, $i=1,\ldots,m$. Put
$$
T'\; =\; R\cup \{\neg (t_1=s_1), \ldots , \neg (t_m=s_m)\}.
$$
Since $T'$ is realized in $\Ce$ and $\Ce\in \ucl(\B)$, then, by
Lemma~\ref{lemmaT0}, $T'$ is finitely satisfiable in $\B$.
However, $T'$ is not satisfiable in $\B$. Thus $\B$ is not
$\ur$-compact. Hence, by Proposition~\ref{localU}, $\Ce$ is not
discriminated by $\B$. We proved $6) \Longrightarrow 4)$.

Equivalence $6) \Longleftrightarrow 7)$ is true in more general
case. Let $\K$ and $\K'$ be two classes of $\L$-algebras (let us
have in mind $\K=\ucl(\B)$ and $\K'=\Dis (\B)$), $\K$ is universal
axiomatizable and $\K'$ is closed under taking $\L$-subalgebras.
Then $\K=\Lo \K'$ is equivalent to $\K_\omega=\K'_\omega$. Indeed,
$\K=\Lo \K'$ easy implies $\K_\omega=\K'_\omega$. Inversely, if
$\K_\omega=\K'_\omega$ then $\K=\Lo \K_\omega=\Lo \K'_\omega=\Lo
\K'$.

\bigskip

Now we begin to prove Theorem~\ref{theoremQ}.

Equivalences $5) \Longleftrightarrow 6)$, $4) \Longleftrightarrow
6)$, $6) \Longleftrightarrow 7)$ may be proven by means of the
similar reasoning as in Theorem~\ref{theoremU} (remind that
$\pvar(\B)=\Res(\B)$).

Let us show equivalence $7) \Longleftrightarrow
8)\Longleftrightarrow 9)$. For an arbitrary algebra $\B$ we have
$\qvar (\B)=\Ls \S \P (\B)=\Ld \S \P
(\B)$~\cite[Corollary~2.3.4]{Gorbunov} and $\Res (\B)=\S\P(\B)$.
So the identity $\qvar(\B)=\Lo\Res (\B)$ is equivalent to $\Ls \S
\P (\B)=\Lo \S \P (\B)$ or $\Ld \S \P (\B)=\Lo \S \P (\B)$.

Equivalence $2) \Longleftrightarrow 10)$ is easy. Equivalence $11)
\Longleftrightarrow 12)$ is due to
V.\,A.\,Gorbunov~\cite[Proposition~1.4.9]{Gorbunov}. So, it
remains to prove implications $2) \Longrightarrow 11)$ and $11)
\Longrightarrow 10)$.

Let $\B$ be $\qw$-compact algebra, $\{S_i, i\in I\}$ a directed
system of radical ideals from $\At_\L(X)$ and $S=\bigcup_{i\in
I}{S_i}$. We show that $S=\Rad (S)$, i.e., $\Rad(S) \subseteq
\bigcup_{i\in I}{S_i}$. Indeed, if $c$ is a consequence of $S$
then there exists a finite subsystem $S_0\subseteq S$ with $c\in
\Rad(S_0)$. Since $I$ is directed there exists an index $i\in I$
such that $S_0 \subseteq S_i$, therefore $c\in S_i$. Thus we have
implication $2) \Longrightarrow 11)$.

To prove implication $11) \Longrightarrow 10)$ consider an
arbitrary system $S\subseteq \At_\L(X)$. The family $\{\Rad
(S_0)\}$, where $S_0$ runs all finite subsystems of a system $S$,
forms a directed system of radical ideals from $\At_\L(X)$. Hence
$\bigcup_{S_0 \subseteq S}{\Rad (S_0)}$ is a radical ideal over
$\B$. Also we have
$$
S \subseteq \bigcup_{S_0 \subseteq S} \: \Rad _\B (S_0) \subseteq
\Rad_\B(S),
$$
therefore $\bigcup_{S_0 \subseteq S}{\Rad (S_0)} = \Rad_\B(S)$.
So, implication $11) \Longrightarrow 10)$ has been proven.

\section{Weakly equationally Noetherian and weakly $\uw$-compact
algebras}\label{sec:weak}

A weak form of the equationally Noetherian property naturally
arises in practice. We discuss algebras with this property in
Subsection~\ref{subsec:weakEN}.

In Subsection~\ref{subsec:weak} we have introduced weakly
$\uw$-compact algebras. Now in Subsection~\ref{subsec:weakcompact}
we present some equivalent approaches to weakly $\uw$-compact
algebras.

In Subsection~\ref{subsec:logirr} we study logically irreducible
algebraic sets. It is important ro note that logically irreducible
algebraic sets inspired the notion of weakly $\uw$-compact
algebras.

\subsection{Weak equationally Noetherian
property}\label{subsec:weakEN}

\begin{definition}
An $\L$-algebra $\B$ is said to be {\em weakly equationally
Noetherian}, if for any finite set $X$ every system $S\subseteq
\At_\L(X)$ is equivalent over $\B$ to some finite system $S_0
\subseteq \At_\L(X)$. Here we do not assume that $S_0$ is a
subsystem of $S$.
\end{definition}

To make comparison equationally Noetherian and weakly equationally
Noetherian properties it is required to reformulate corresponding
definitions in the following form.

An $\L$-algebra $\B$ is termed {\em weakly equationally
Noetherian}, if for any finite set $X$ and any system $S\subseteq
\At_\L(X)$ there exists finite system $S_0 \subseteq \Rad_\B(S)$
such that $\V_\B(S)=\V_\B(S_0)$.

An $\L$-algebra $\B$ is termed {\em equationally Noetherian}, if
for any finite set $X$ and any system $S\subseteq \At_\L(X)$ there
exists finite system $S_0 \subseteq [S]$ such that
$\V_\B(S)=\V_\B(S_0)$.

Indeed, for every atomic formula $c=(t=s)\in [S]$ there exists a
finite subsystem $S_c\subseteq S$ such that $S_c \vdash (t=s)$.
Therefore, if $\V_\B(S)=\V_\B(S_0)$ for a finite system $S_0
\subseteq [S]$ then one has
\begin{equation}\label{eq:V}
 \V_\B(S)=\V_\B(\bigcup_{c\in S_0}
{S_c}).
\end{equation}

\begin{lemma}\label{lemma4}
If an $\L$-algebra $\B$ is weakly equationally Noetherian and
$\qw$-compact then it is equationally Noetherian.
\end{lemma}

\begin{proof}
As $\B$ is weakly equationally Noetherian, for each system of
equations $S$ there exists a finite system $S_0 \subseteq
\Rad_\B(S)$ with $\V_\B(S)=\V_\B(S_0)$. As $\B$ is $\qw$-compact,
for each equation $c=(t_0=s_0)\in S_0$ there exists a finite
subsystem $S_c\subseteq S$ with $\V_\B(S_c) \subseteq
\V_\B(t_0=s_0)$. Thereby, one has~\eqref{eq:V}. It means that $\B$
is equationally Noetherian algebra.
\end{proof}

\begin{lemma}
If an $\L$-algebra $\B$ is weakly equationally Noetherian and
$\Ce$ a subalgebra of some  direct power of $\Ce$ then $\Ce$ is
weakly equationally Noetherian too.
\end{lemma}

\begin{proof}
It follows from~\cite[Lemma~3.7]{DMR2}.
\end{proof}

It is clear that every weakly equationally Noetherian algebra is
$\E$-compact.

\begin{lemma}\label{lemma3}
If an $\L$-algebra $\B$ is weakly equationally Noetherian then
$$
\ucl(\B)\; \cap\; \Res(\B)_\omega = \Dis(\B)_\omega.
$$
\end{lemma}

\begin{proof}
Since $\Dis(\B)\subseteq\Res(\B)$, $\Dis(\B)\subseteq\ucl(\B)$,
and $\Res(\B)=\pvar(\B)$ for any algebra $\B$~\cite{DMR1}, we
should check that $\ucl(\B)\, \cap\, \pvar(\B)_\omega \subseteq
\Dis(\B)_\omega$. Let us assume that $\Ce$ is a finitely generated
algebra such that $\Ce\in\pvar(\B)\setminus\Dis(\B)$ and prove
$\Ce\not\in \ucl(\B)$. If $\Ce$ is the trivial algebra $\E$ then,
by definition, condition $\Ce\not\in\Dis(\B)$ implies that $\B$
has not a trivial subalgebra. Since $\B$ is weakly equationally
Noetherian, then $\B$ is $\E$-compact, and, by Lemma~\ref{lemma2},
$\E\not\in\ucl(\B)$. Thereby, we may assume that $\Ce$ is
non-trivial.

Let $\langle \{c_1,\ldots,c_n\} \mid S \rangle$ be a presentation
of $\Ce$, i.e., $\Ce \simeq \T_\L(X)/\theta_S$, $S\subseteq \At_\L
(X)$, $X=\{x_1,\ldots,x_n\}$. Since $\Ce\not\in\Dis(\B)$, there
exits atomic formulas $(t_1=s_1),\ldots,(t_m=s_m)\in
\At_\L(X)\setminus [S]$ such that the (infinite) formula
$$
\forall \;  x_1 \ldots \forall \;  x_n \left( \bigwedge
\limits_{(t=s) \in S} t (\bar{x})=s (\bar{x}) \; \longrightarrow
\; \bigvee\limits_{i=1}^m t_i (\bar{x})=s_i(\bar{x}) \right)
$$
holds in $\B$. As one can find a finite system $S_0\subseteq
\At_\L(X)$ with $\V_\B(S_0)=\V_\B(S)$ then the universal sentence
\begin{equation}\label{eq:S0}
\forall \;  x_1 \ldots \forall \;  x_n \left(\bigwedge
\limits_{(t=s) \in S_0} t (\bar{x})=s (\bar{x}) \; \longrightarrow
\; \bigvee\limits_{i=1}^m t_i (\bar{x})=s_i(\bar{x}) \right)
\end{equation}
holds in $\B$.

Since $\V_\B(S_0)=\V_\B(S)$ we have
$\V_\Ce(S_0)=\V_\Ce(S)$~\cite[Lemma~3.7]{DMR2}. Hence,
$(c_1,\ldots,c_n)\in \V_\Ce(S_0)$ but $t_i(c_1,\ldots,c_n)\ne
s_i(c_1,\ldots,c_n)$ for all $i=1,\ldots,m$. Therefore, universal
formula~\eqref{eq:S0} is not true in $\Ce$, and
$\Ce\not\in\ucl(\B)$.
\end{proof}

\subsection{Logically irreducible algebraic
sets}\label{subsec:logirr}

One of the approaches to $\uw$-compact algebras deals with
so-called logically irreducible algebraic sets.

\begin{definition}
We say that an algebraic set $Y$ over $\B$ is {\em logically
irreducible} if its coordinate algebra $\Gamma (Y)$ belongs to
$\ucl(\B)$.
\end{definition}

In Section~\ref{sec:unification_theorems_compact} we have
discussed that every irreducible algebraic set over an arbitrary
algebra  $\B$ is logically irreducible. In
Subsection~\ref{subsec:weakcompact} we will show that the notions
of irreducible and logically irreducible algebraic sets coincide
if and only if $\B$ is weakly $\uw$-compact algebra.

\begin{lemma}
Let $\B$ be an $\L$-algebra. For a finitely generated $\L$-algebra
$\Ce$ the following conditions are equivalent:
\begin{itemize}
\item $\Ce$ is the coordinate algebra of a logically irreducible algebraic set over $\B$;
\item $\Ce$ belongs to $\ucl(\B)\, \cap \, \pvar(\B)$.
\end{itemize}
\end{lemma}

\begin{proof}
Indeed, $\Ce$ is the coordinate algebra of an algebraic set over
$\B$ if and only if $\Ce\in
\pvar(\B)$~\cite[Proposition~3.22]{DMR2}.
\end{proof}

\begin{corollary}\label{lemma1}
The class of all coordinate algebras of logically irreducible
algebraic sets over $\B$ coincides with $\ucl(\B)\: \cap \:
\pvar(\B)_\omega$.
\end{corollary}

For irreducible algebraic sets we have the following result.

\begin{lemma}[\cite{DMR2}]\label{lemma3.42}
Let $\B$ be an $\L$-algebra. Every non-empty algebraic set $Y$
over $\B$ is a union of maximal with respect to inclusion
irreducible algebraic sets over $\B$.
\end{lemma}

Now we try to find a similar decomposition for algebraic sets into
a union of maximal logically irreducible algebraic sets. It is
clear that Lemma~\ref{lemma3.42} gives a decomposition. However,
maximal with respect to inclusion irreducible algebraic set may be
a proper subset of some logically irreducible algebraic set.

\begin{lemma}\label{irr9}
Let $Y_1\subset Y_2 \subset \ldots$ be an ascending chain of
logically irreducible algebraic sets in $B^n$ and $Y$ the least
algebraic set containing all these sets. Then $Y$ is logically
irreducible algebraic set.
\end{lemma}

\begin{proof}
Note that  $Y=\V_\B (\Rad (\bigcup\limits_{i}{Y_i}))$ and $\Rad
(Y)=\bigcap\limits_{i}{\Rad (Y_i)}$. Hence, there exists embedding
$h\colon \Gamma (Y) \to \prod\limits_{i}{\Gamma
(Y_i)}$~\cite[Lemma~3.1]{DMR1}. Index $i$ runs the linearly
ordered set $I$. For each $i\in I$ denote by $J_i$ the set $\{j\in
I, j\geqslant i\}$. The family of subsets $\{J_i, i\in I\}$ is
centered, hence there exists an ultrafilter $D$ on $I$ containing
$J_i$ for all $i\in I$. Let $f\colon \prod\limits_{i}{\Gamma
(Y_i)}\to \prod\limits_{i}{\Gamma (Y_i)}/D$ be a canonical
homomorphism. Let us show that composition $f\circ h\colon \Gamma
(Y) \to \prod\limits_{i}{\Gamma (Y_i)}/D$ is embedding.

Indeed, we have $\Gamma (Y)=\T _\L(X)/\theta_{\Rad (Y)}$, where
$X=\{x_1,\ldots,x_n\}$. If $t_1/\theta_{\Rad (Y)},
t_2/\theta_{\Rad (Y)}$ are distinct elements from $\Gamma (Y)$
then $(t_1=t_2)\in \At _\L (X)\setminus \Rad(Y)$. Since $\Rad
(Y_1) \supset \Rad (Y_2) \supset \ldots$, then there exists an
index $i_0\in I$ such that $(t_1=t_2)\not\in\Rad(Y_i)$ for all
$i\in J_{i_0}$. It implies that $f\circ h (t_1/\theta_{\Rad (Y)})
\ne f\circ h (t_2/\theta_{\Rad (Y)})$. Thus $f\circ h$ is
injective.

Since $\Gamma (Y_i) \in \ucl (\B)$ for each $i\in I$ and $\Gamma
(Y) \in \S \Pu (\{\Gamma (Y_i), i\in I\})$, then $\Gamma (Y) \in
\ucl (\B)$, i.e., $Y$ is logically irreducible algebraic set.
\end{proof}

\begin{lemma}
Let $\B$ be an $\L$-algebra.  Every non-empty algebraic set $Y$
over $\B$ is a union of maximal with respect to inclusion
logically irreducible algebraic sets over $\B$.
\end{lemma}

\begin{proof}
We will show that for each point $p\in Y$ there exists logically
irreducible algebraic set $Z$ such that $p\in Z \subseteq Y$ and
$Z$ is maximal with these properties. Denote by $\Omega$ the
family of logically irreducible algebraic sets $Z$ with $p\in Z
\subseteq Y$ and show that $\Omega$ is not empty and has maximal
elements.

Denote by $Z_p$ the closure in the Zariski topology of the set
$\{p\}$. One has $p\in Z_p\subseteq Y$. Furthermore, $Z_p$ is
irreducible algebraic set~\cite[Lemma~3.34]{DMR2}. Hence, $Z_p\in
\Omega$.

By Zorn Lemma it is sufficiently to show now that family $\Omega$
contains upper boundary for each ascending chain $Y_1\subset Y_2
\subset \ldots$ of element from $\Omega$. Let $Y_p$ be the least
algebraic set that contains union $\bigcup\limits_{i}{Y_i}$. By
Lemma~\ref{irr9}, $Y_p$ is logically irreducible. As $Y_p\subseteq
Y$ one has $Y_p \in \Omega$.

Thereby, the union $\bigcup_{p\in Y} {Y_p}$ is desired.
\end{proof}

Let us remind that for equationally Noetherian algebras we have
the next result.

\begin{theorem}[\cite{DMR1}]\label{irr}
Let $\B$ be an equationally Noetherian algebra. Then any non-empty
algebraic set $Y$ over $\B$ is a finite union of irreducible
algebraic sets (irreducible components): $Y=Y_1 \cup \ldots \cup
Y_m$. Moreover, if  $Y_i \not \subseteq Y_j$ for $i\ne j$ then
this decomposition is unique up to a permutation of components.
\end{theorem}

It is natural to ask the following question.

\noindent{\bf Decomposition Problem.} {\it Let $\B$ be a ``good''
algebra ($\uw$-,$\qw$-compact, weakly equationally Noetherian, for
instance). Is it true that every non-empty algebraic set over $\B$
is a finite union of logically irreducible algebraic sets?}

In spite of the fact that $\uw$-compact and weakly equationally
Noetherian algebras are the closest algebras to equationally
Noetherian ones we give for them the negative answer to the
question above.

Indeed, a decomposition $Y=Y_1 \cup \ldots \cup Y_m$ of algebraic
set $Y$ into a union of algebraic sets $Y_1,\ldots,Y_m$ implies
the existence of a subdirect embedding $h\colon \Gamma(Y)\to
\Gamma(Y_1)\times\ldots\times\Gamma(Y_m)$~\cite{DMR2}. Suppose
that the Decomposition Problem has the positive answer for
$\uw$-compact algebras. It involves that the Embedding Problem for
$\uw$-compact algebras has the positive answer too. However,
A.\,N.\,Shevlyakov has proven the inverse result (see
Subsection~\ref{subsec:C}). Moreover, he has proven also that the
Decomposition Problem for weakly equationally Noetherian algebras
has the negative answer~\cite{Shevl1}.

\subsection{Weak $\uw$-compactness}\label{subsec:weakcompact}

In the proposition below we gather the different approaches to
weakly $\uw$-compact algebras.

\begin{proposition}\label{propwU}
For an algebra $\B$ in a functional language $\L$ the following
conditions are equivalent:
\begin{itemize}
\item[1)] $\B$ is weakly $\uw$-compact;
\item[2)] every non-empty logically irreducible algebraic set over $\B$ is
irreducible;
\item[3)] every non-trivial coordinate algebra over $\B$ that
belongs to $\ucl(\B)$ is irreducible;
\item[4)] $\ucl(\B)\,\cap\,\Res(\B)_\omega=(\Dis(\B)\e)_\omega$.
\end{itemize}
\end{proposition}

\begin{proof}
Equivalence $1) \Longleftrightarrow 2)$ is evident by definition.
Remind that the trivial algebra $\E$ is a coordinate algebra over
$\B$ anyway, moreover, if $Y$ is an algebraic set over $\B$ such
that $\E=\Gamma (Y)$ then $Y$ is irreducible or
$Y=\emptyset$~\cite[Lemma~3.22]{DMR2}. It implies that we have
equivalence $2) \Longleftrightarrow 3)$.

Since $\Dis(\B)\subseteq\Res(\B)$, $\Dis(\B)\subseteq\ucl(\B)$,
$\Res(\B)=\pvar(\B)$, and $\E \in \Res(\B)$ for any algebra $\B$,
then item~4) means that every non-trivial algebra $\Ce$ from $\ucl
(\B) \cap \pvar (\B)_\omega$ belongs to $\Dis(\A)_\omega$.

As the class of all coordinate algebras of irreducible algebraic
sets over $\B$ coincides with
$\Dis(\B)_\omega$~\cite[Corollary~3.37]{DMR2}, and, by
Corollary~\ref{lemma1}, the class of all coordinate algebras of
logically irreducible algebraic sets over $\B$ coincides with
$\ucl(\B)\: \cap \: \pvar(\B)_\omega$, we have equivalence $3)
\Longleftrightarrow 4)$.
\end{proof}

\begin{remark}
Every $\uw$-compact (as well as $\qw$-compact, weakly equationally
Noetherian) algebra is $\E$-compact. However, there exist weakly
$\uw$-compact algebras that are not $\E$-compact (see
Example~\ref{example} bellow). Suppose an algebra $\B$ is
$\E$-compact. In this case one can omit ``non-empty'' in item~2),
omit ``non-trivial'' in item~3), and write
``$\ucl(\B)\,\cap\,\Res(\B)_\omega=\Dis(\B)_\omega$'' instead of
``$\ucl(\B)\,\cap\,\Res(\B)_\omega=(\Dis(\B)\e)_\omega$'' in
item~4) in the formulation of Proposition~\ref{propwU}. In this
case the empty set is not algebraic over $\B$, or if it is
algebraic then its coordinate algebra $\E$ does not belong to
$\ucl (\B)$.
\end{remark}

\begin{lemma}\label{lemma5}
If an $\L$-algebra $\B$ is weakly $\uw$-compact and $\qw$-compact
then it is $\uw$-compact.
\end{lemma}

\begin{proof}
We need to show that $\ucl (\B)_\omega \subseteq \Dis(\B)_\omega$.
Assume that $\Ce$ is a finitely generated algebra and
$\Ce\not\in\Dis(\B)$. Since $\Dis(\B)_\omega=\ucl(\B)\, \cap\,
\pvar(\B)_\omega$, then $\Ce\not \in\ucl(\B)$, and we have
required, or $\Ce\not\in\pvar(\B)$. By Theorem~\ref{theoremQ},
$\Ce\not\in\pvar(\B)$ implies that $\Ce\not\in\qvar(\B)$, hence
$\Ce\not \in\ucl(\B)$.
\end{proof}

The next question is naturally arises. Is there a geometric
definition of weak $\uw$-compactness?

\begin{definition}
We name an $\L$-algebra $\B$ {\em geometrically weakly
$\uw$-compact} if for any finite set $X$, any system of equations
$S\subseteq \At_\L (X)$, and any equations $(t_1=s_1), \ldots,
(t_m=s_m)\in \At_\L (X)$ such that
$$
\V_\B(S) \; \subseteq \; \V_\B(t_1=s_1)\:\cup \:\ldots\: \cup\:
\V_\B(t_m=s_m)
$$
and for each $i\in\{1,\ldots,m\}$
$$
\V_\B(S) \; \nsubseteq \; \V_\B(t_i=s_i)
$$
there exists a finite subsystem $S_0\subseteq \Rad_\B(S)$ such
that
$$
\V_\B(S_0) \; \subseteq \; \V_\B(t_1=s_1)\:\cup \:\ldots\: \cup\:
\V_\B(t_m=s_m).
$$
\end{definition}

The definition above is evident generalization of both weak
equationally Noetherian property and $\uw$-compactness. It also
has analogs in terms of radical, in terms of infinite formulas,
and in terms of compactness.

\begin{lemma}\label{theoremwU}
For an algebra $\B$ in a functional language $\L$ the following
conditions are equivalent:
\begin{itemize}
\item[1)] $\B$ is geometrically weakly $\uw$-compact;
\item[2)] for any finite set $X$, any radical ideal $S\subseteq \At_\L
(X)$ over $\B$, and any atomic formulas $(t_1=s_1), \ldots,
(t_m=s_m)\in \At_\L (X)\setminus \Rad_\B (S)$ if
$c=(t_1=s_1)\vee\ldots \vee (t_m=s_m)$ is a consequence of $S$
over $\B$ then there exists a finite subsystem $S_c\subseteq S$
such that $c$ is a consequence of $S_c$ over $\B$;
\item[3)] for any finite set $X$, any radical ideal $S\subseteq \At_\L
(X)$ over $\B$, and any atomic formulas $(t_1=s_1), \ldots,
(t_m=s_m)\in \At_\L (X)$ if an (infinite) formula
$$
\forall \;  x_1 \ldots \forall \;  x_n \left( \bigwedge
\limits_{(t=s) \in S} t (\bar{x})=s (\bar{x}) \; \longrightarrow
\; \bigvee\limits_{i=1}^m t_i (\bar{x})=s_i(\bar{x}) \right)
$$
holds in $\B$, and for each $i \in\{1,\ldots,m\}$ an (infinite)
formula
$$
\forall \;  x_1 \ldots \forall \;  x_n \left( \bigwedge
\limits_{(t=s) \in S} t (\bar{x})=s (\bar{x}) \; \longrightarrow
\; t_i(\bar{x})=s_i(\bar{x}) \right)
$$ does not hold in $\B$, then for some finite subsystem $S_c\subseteq S$ the universal
sentence
$$
\forall \;  x_1 \ldots \forall \;  x_n \left( \bigwedge
\limits_{(t=s) \in S_c} t (\bar{x})=s (\bar{x}) \; \longrightarrow
\; \bigvee\limits_{i=1}^m t_i (\bar{x})=s_i(\bar{x}) \right)
$$
holds in $\B$;
\item[4)] for any finite set $X$, any radical ideal $S\subseteq \At_\L
(X)$ over $\B$, and any atomic formulas $(t_1=s_1), \ldots,
(t_m=s_m)\in \At_\L (X)$ if the set of formulas
$$
T\; =\; S\cup \{\neg (t_1=s_1), \ldots , \neg (t_m=s_m)\}
$$
is finitely satisfiable in $\B$ and for each $i \in\{1,\ldots,m\}$
the set of formulas
$$
T_i\; =\; S\cup \{\neg (t_i=s_i)\}
$$
is realized in $\B$ then $T$ is satisfiable in $\B$.
\end{itemize}
\end{lemma}

\begin{proof}
Equivalences $1) \Longleftrightarrow 2)$, $1) \Longleftrightarrow
3)$, $3) \Longleftrightarrow 4)$ follows from
Lemma~\ref{lemmacompact}. Note that the statement in item~3) has a
form ``$A \, \& \, \neg C$ implies $B$''. The equivalent statement
is ``$\neg B\, \& \, \neg C$ implies $\neg A$'' which gives~4). So
we have $3) \Longleftrightarrow 4)$.
\end{proof}

Unfortunately, for weak $\uw$-compactness we have no an analog of
Theorem~\ref{theoremU} that holds for $\uw$-compact algebras.

\begin{lemma}
If an $\L$-algebra $\B$ is geometrically weakly $\uw$-compact that
it is weakly $\uw$-compact. The converse statement does not hold.
\end{lemma}

\begin{proof}
Suppose that $\B$ is geometrically weakly $\uw$-compact and $Y$ a
non-empty algebraic set over $\B$ such that $\Gamma(Y)\in \ucl
(\B)$. We need to show that $\Gamma (Y) \in \Dis (\B)$. Let
$S=\Rad(Y)$, then $\Gamma (Y)$ has the presentation $\langle X
\mid S \rangle$. If $\Gamma (Y)$ is the trivial algebra, i.e.,
$S=\At_\L(X)$, then $Y$ is irreducible~\cite[Lemma~3.22]{DMR2} and
$\Gamma (Y) \in \Dis (\B)$.

Assume now that $\Gamma (Y)$ is non-trivial, i.e.,
$S\ne\At_\L(X)$. As the coordinate algebra $\Gamma (Y)$ is
separated by $\B$, hence for each atomic formula $(t=s)\in
\At_\L(X)\setminus S$ the set of formulas $S\cup \{\neg (t=s)\}$
is realized in $\B$. Take atomic formulas $(t_1=s_1), \ldots,
(t_m=s_m)\in \At_\L (X)\setminus S$. As the set of formulas $T
=S\cup \{\neg (t_1=s_1), \ldots , \neg (t_m=s_m)\}$ is satisfiable
in $\langle X \mid S \rangle$, and $\langle X \mid S
\rangle\in\ucl(\B)$, then, by Lemma~\ref{lemmaT0}, $T$ is finitely
satisfiable in $\B$. It follows from item~4) of
Lemma~\ref{theoremwU} that $T$ is satisfiable in $\B$. Thereby,
algebra $\langle X \mid S \rangle$ is discriminated by $\B$.

Example~\ref{example} below shows that the converse statement does
not hold.
\end{proof}

The following example is similar to the example by
M.\,V.\,Kotov~\cite{Kotov}.

\begin{example}\label{example}
Let $\L=\{g_n, n\in \mathbb{N}\}$ be the infinite signature with
unary functional symbols and $\A$ the $\L$-algebra with the
universe $\mathbb{N}$ and
$$
g_n(x) = \left\{
\begin{array}{ll}
2n,     & x = 2n + 1, \\
2n + 1, & x = 2n, \\
x, & \mbox{otherwise}.
\end{array}
\right.
$$
It is clear that $\A$ has no trivial subalgebra. At the same time,
the set of formulas $\{g_n(x)=x, n \in \mathbb{N}\}$ is finitely
satisfiable in $\A$, therefore, by Compactness Theorem, it is
satisfiable in some ultrapower $\A^\ast$ of $\A$. As $\A^\ast \in
\ucl (\A)$, then $\E \in \ucl (\A)$. Thereby, $\A$ is not
$\E$-compact.

We state that $\A$ is weakly $\uw$-compact. Indeed, take a
non-trivial algebra $\Ce$ from $\ucl (\A) \cap \pvar (\A)$. Since
$\pvar(\A)=\S\P(\A)$ then $\Ce$ is a subalgebra of a direct power
of $\A$. For any $n,m\in \mathbb{N}$, $n\ne m$, the universal
formula
$$
\forall \; x  \quad (\, g_n(x)=x \; \vee \; g_m(x)=x \, )
$$
holds in $\A$. Therefore, $\Ce$ has a finite universe $\{c_1,c'_1,
\ldots, c_d, c'_d\}$ with $c_i=g_{n_i}(c'_i)$ for all
$i=\overline{1,d}$. The map $h\colon \Ce \to \A$, $h(c_i)=2n_i$,
$h(c'_i)=2n_i+1$, $i=\overline{1,d}$, is a monomorphism. Thus,
$\Ce \in \Dis(\A)$, and $\A$ is weakly $\uw$-compact.

Let us check that $\A$ is not geometrically weakly $\uw$-compact.
Consider the systems of equations $S'(x)=\{g_n(x)=x,\: n\in
\mathbb{N}\setminus\{0\}\}$ and $S(x,y)=S'(x)\cup S'(y)$. We have
$\V_\A(S)=\{(0,0), (0,1), (1,0), (1,1)\}$. Therefore,
\begin{gather*}
\V_\A(S) \; \subseteq \; \V_\A(x=y) \: \cup \: \V_\A(x=g_0(y)),\\
\V_\A(S)\; \nsubseteq \; \V_\A(x=y), \quad \V_\A(S) \; \nsubseteq
\; \V_\A(x=g_0(y)).
\end{gather*}
Furthermore, it is not hard to see that
$$
\Rad_\A(S)= \left\{
\begin{array}{l}
x=g_{n_1}(g_{n_2}(\ldots g_{n_m}(x)\ldots)),  \\
y=g_{n_1}(g_{n_2}(\ldots g_{n_m}(y)\ldots)), \; n_i\ne 0\\
x=g_{n_1}(g_{n_2}(\ldots g_{n_m}(y)\ldots)).
\end{array}
\right\}.
$$
It is obvious that for any finite subsystem $S_0 \subseteq
\Rad_\A(S)$ we have
$$
\V_\A(S_0) \; \nsubseteq \; \V_\A(x=y) \: \cup \: \V_\A(x=g_0(y)).
$$
\end{example}

\section{Connections between the classes of algebras $\QC$, $\UC$, $\UC'$, $\EN'$, and $\EN$}\label{sec:picture}

Let $\L$ be a functional language. We use the following
denotations:\\
$\EN\,\:$~--- the class of all equationally
Noetherian $\L$-algebras;\\
$\EN\,'$~--- the class of all weakly
equationally Noetherian $\L$-algebras;\\
$\QC\,\:$~--- the class of all $\qw$-compact $\L$-algebras;\\
$\UC\,\:$~--- the class of all $\uw$-compact $\L$-algebras;\\
$\UC'$~--- the class of all weakly $\uw$-compact $\L$-algebras.

It is clear that
$$
\QC\;  \supseteq \; \UC \; \supseteq \; \EN \; \subseteq \; \EN'.
$$
Moreover, by Lemma~\ref{lemma4},
$$
\EN  =\EN'\cap \QC = \EN'\cap \UC.
$$

So, we have exactly the following picture for co-location of
classes $\EN$, $\EN'$, $\QC$, $\UC$:
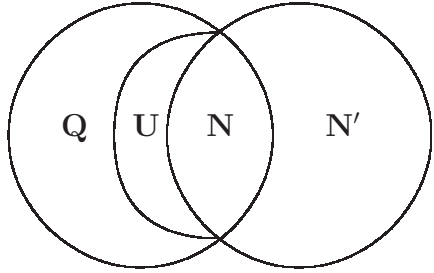
\begin{center}
\begin{picture}(160,100)
\qbezier[1100](50,0)(69,0)(85,15)
\qbezier[1100](85,15)(100,31)(100,50)
\qbezier[1100](100,50)(100,69)(85,85)
\qbezier[1100](85,85)(69,100)(50,100)
\qbezier[1100](50,100)(31,100)(15,85)
\qbezier[1100](15,85)(0,69)(0,50)
\qbezier[1100](0,50)(0,31)(15,15)
\qbezier[1100](15,15)(31,0)(50,0)

\qbezier[1100](110,0)(129,0)(145,15)
\qbezier[1100](145,15)(160,31)(160,50)
\qbezier[1100](160,50)(160,69)(145,85)
\qbezier[1100](145,85)(129,100)(110,100)
\qbezier[1100](110,100)(91,100)(75,85)
\qbezier[1100](75,85)(60,69)(60,50)
\qbezier[1100](60,50)(60,31)(75,15)
\qbezier[1100](75,15)(91,0)(110,0)

\qbezier[1100](40,50)(40,89)(80,89)
\qbezier[1100](40,50)(40,11)(80,11)

\put(75,50){$\EN$} \put(120,50){$\EN'$} \put(20,50){$\QC$}
\put(47,50){$\UC$}
\end{picture}
\end{center}


Let us find the place of the class $\UC'$ in the picture above. By
Theorem~\ref{uw}, Lemma~\ref{lemma3} and Proposition~\ref{propwU},
we have
$$
\UC\subseteq  \UC' \quad \mbox{and} \quad \EN'\subseteq  \UC'.
$$
It follows from Lemma~\ref{lemma5} that
$$
\QC\, \cap\, \UC'=\UC.
$$

Hence, co-location of classes $\EN$, $\EN'$, $\QC$, $\UC$, and
$\UC'$ are exactly the following:
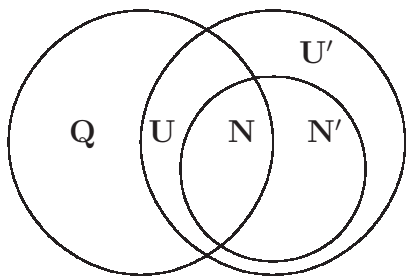
\begin{center}
\begin{picture}(150,100)
\qbezier[1100](50,0)(69,0)(85,15)
\qbezier[1100](85,15)(100,31)(100,50)
\qbezier[1100](100,50)(100,69)(85,85)
\qbezier[1100](85,85)(69,100)(50,100)
\qbezier[1100](50,100)(31,100)(15,85)
\qbezier[1100](15,85)(0,69)(0,50)
\qbezier[1100](0,50)(0,31)(15,15)
\qbezier[1100](15,15)(31,0)(50,0)

\qbezier[1100](100,0)(119,0)(135,15)
\qbezier[1100](135,15)(150,31)(150,50)
\qbezier[1100](150,50)(150,69)(135,85)
\qbezier[1100](135,85)(119,100)(100,100)
\qbezier[1100](100,100)(81,100)(65,85)
\qbezier[1100](65,85)(50,69)(50,50)
\qbezier[1100](50,50)(50,31)(65,15)
\qbezier[1100](65,15)(81,0)(100,0)

\qbezier[1100](100,5)(114,5)(125,15)
\qbezier[1100](125,15)(135,26)(135,40)
\qbezier[1100](135,40)(135,54)(125,65)
\qbezier[1100](125,65)(114,75)(100,75)
\qbezier[1100](100,75)(86,75)(75,65)
\qbezier[1100](75,65)(65,54)(65,40)
\qbezier[1100](65,40)(65,26)(75,15)
\qbezier[1100](75,15)(86,5)(100,5)

\put(23,50){$\QC$} \put(53,50){$\UC$}\put(83,50){$\EN$}
\put(113,50){$\EN'$}\put(110,80){$\UC'$}
\end{picture}
\end{center}

In paper~\cite{MR2} A.\,G.\,Myasnikov and V.\,N.\,Remeslennikov
asked the questions for the class of groups:

{\bf Question 1:} \;\quad\quad { $\EN=\QC$ \quad \quad or
\quad\quad  $\EN \ne \QC$ ?}

{\bf Question 2:} \;\quad\quad { $\QC=\UC$ \quad \quad or
\quad\quad $\QC \ne \UC$ ?}

Now we add new questions:

{\bf Question 3:} \quad\quad { $\EN=\EN'$ \quad \quad or
\quad\quad $\EN \ne \EN'$ ?}

{\bf Question 4:} \quad\quad { $\EN=\UC$ \:\quad \quad or
\quad\quad $\EN \ne \UC$ ?}

{\bf Question 5:} \; { $\UC'=\UC \:\cup\: \EN'$ \; or \;
$\UC'\ne\UC\: \cup\: \EN'$ ?}

The answer to the first question has been given by B.\,I.\,Plotkin
in~\cite{Plot3}. He has constructed $\qw$-compact group that is
not equationally Noetherian. We will discuss that construction in
this section below. Note that B.\,I.\,Plotkin uses notation {\em
logically Noetherian} for $\qw$-compact algebras and {\em
geometrically Noetherian} for equationally Noetherian algebras.

The second and third questions have been solved by
M.\,V.\,Kotov~\cite{Kotov}. He has constructed examples that show
$\QC \ne \UC$ and $\EN\ne\EN'$. His examples are original
algebraic structures in the language $\L=\{g_n, n \in
\mathbb{N}\}$ with countable set of unary functional symbols and
with universe-sets $\mathbb{R}$ and $\mathbb{N}$.

At these results the fourth question remains open as well as the
problem of differentiation of classes $\QC$, $\UC$, $\EN$, $\EN'$
for classical varieties: groups, rings, monoids, semigroups.
In~\cite{Shevl1} A.\,N.\,Shevlyakov finds the neat examples in the
variety of commutative idempotent semigroups in the language with
countable set of constants. His examples distinguish classes
$\EN$, $\EN'$, $\QC$, $\UC$.

The algebra $\A$ from Example~\ref{example} gives an answer to the
fifth question. It has been shown that $\A \in \UC'$, but $\A$ is
not $\E$-compact. Since all $\uw$-compact and weakly equationally
Noetherian algebras are $\E$-compact, then $\A \not \in
\UC\,\cup\,\EN'$. Another example for $\UC'\ne\UC\: \cup\: \EN'$
has been constructed by A.\,N.\,Shevlyakov~\cite{Shevl1} in the
class of commutative idempotent semigroups in the language $\L$
with countable set of constants. It is important to note that all
algebras in the language $\L$ are $\E$-compact, by
Lemma~\ref{lemmaL}.

Let us return to the construction given by B.\,I.\,Plotkin. He
denotes by $H$ the discrete direct product of all finitely
generated groups (in the language of groups $\L=\{\cdot, ^{-1},
e\}$). Since every finitely generated group $G$ imbeds into $H$,
then $G$ is a coordinate group over $H$. By~10) in
Theorem~\ref{theoremQ} below, $H$ is $\qw$-compact. As there
exists a finitely generated group $G$ that is not finitely
presented, hence $H$ is not equationally Noetherian.

It is evident that this construction of $H$ may be repeated in
other varieties of algebras, where exist finitely generated, not
finitely presented algebras. Clearly, the algebraic geometry over
objects like $H$ is quite elementary.

\section{$\qw$-compact and $\uw$-compact extensions}
\label{sec:extension}

In Introduction it is given the formulation of the Compactness
Theorem and the notion of logical compactness. The Compactness
Theorem has a great importance in model theory~\cite{Hodges}.

For an arbitrary algebra $\B$ it is possible with a use of the
Compactness Theorem to construct an elementary extension $\B^\ast$
of $\B$ such that $\B^\ast$ is logically compact. This algorithm
is close to the building of the algebraic closure to a given field
$k$.

We use this idea to construct $\uw$-compact elementary extension
for an arbitrary algebra $\B$. At first, let us remind some more
facts from model theory.

\begin{theorem}[Corollary from Los'
Theorem~\cite{Gorbunov}]\label{elem_emb} If $\B^I/D$ is an
ultrapower of an algebra $\B$ then the diagonal map $d\colon \B
\to \B^I/D$, where $d(x)=\bar{x}/D$ and $\bar{x}(i)=x$ for all $i
\in I$, is an elementary embedding.
\end{theorem}

\begin{proposition}[\cite{Marker}]\label{elem_chain} Suppose that $(I, <)$ is a
linear order and $(\M_i, i\in I)$ is an elementary chain. Then
$\M=\bigcup_{i\in I}{\M_i}$ is an elementary extension of each
$\M_i$.
\end{proposition}

Denote by $\mathbb{T}$ the family of all sets of formulas
$$
T\; =\; S\cup \{\neg (t_1=s_1), \ldots , \neg (t_m=s_m)\},
$$
where $S\subseteq \At_\L (X)$, $(t_1=s_1), \ldots, (t_m=s_m)\in
\At_\L (X)$, $\vert X \vert < \infty$. For a given $\L$-algebra
$\B$ let $\mathbb{T}(\B)$ be such subfamily of $\mathbb{T}$ that
$T\in \mathbb{T}(\B)$ if and only if $T$ is finitely satisfiable
in $\B$ but not realized in $\B$. So, algebra $\B$ is
$\uw$-compact if and only if $\mathbb{T}(\B)=\emptyset$.

For $\L$-algebras $\B$ and $\Ce$ we write $\B \equiv_\forall \Ce$
if $\B$ and $\Ce$ are universally equivalent, i.e., $\ucl
(\B)=\ucl(\Ce)$.

\begin{lemma}\label{lemmaT}
Let $\B$ and $\Ce$ be $\L$-algebras and $\B \leq \Ce$. If $\B
\equiv_\forall \Ce$ then $\mathbb{T}(\Ce)\subseteq
\mathbb{T}(\B)$.
\end{lemma}

\begin{proof}
Suppose $T\in \mathbb{T}$ and $T$ is finitely satisfiable in
$\Ce$. Then, by Lemma~\ref{lemmaT0}, $T$ is finitely satisfiable
in $\B$. If $T \not\in \mathbb{T}(\B)$, then $T$ is realized in
$\B$. As $\B \leq \Ce$, then $T$ is realized in $\Ce$ and $T
\not\in \mathbb{T}(\Ce)$.
\end{proof}

\begin{theorem}\label{el_uw}
Let $\B$ be an $\L$-algebra. Then there exists an elementary
extension $\B^\ast$ of $\B$, such that $\B^\ast$ is $\uw$-compact
(in particularly, $\B^\ast$ is $\qw$-compact).
\end{theorem}

\begin{proof}
Consider a well-ordering $(I, <)$ on $\mathbb{T}(\B)$. Let us
construct an elementary chain $(\B_i, i\in I)$. At first, take
$\B_0=\B$. Then $\B_1$ is an ultrapower of $\B$ where $T_0$ is
realized. By Compactness Theorem, such $\B_1$ exists and, by
Theorem~\ref{elem_emb}, $\B_1$ is an elementary extension of $\B$.
Further, $\B_2$ is an ultrapower of $\B_1$ where $T_1$ is
realized, and so on. For an ordinal $\alpha=\beta+1$ we put
$\B_\alpha$ as an ultrapower of $\B_\beta$ where $T_\beta$ is
realized, and $\B_\alpha=\bigcup_{\beta < \alpha}{\B_\beta}$ for a
limit ordinal $\alpha$. Desired algebra $\B^\ast$ is
$\bigcup_{i\in I}{\B_i}$. Indeed, $\B^\ast$ is an elementary
extension of $\B$, by Theorem~\ref{elem_emb} and
Proposition~\ref{elem_chain}.

Let us show that $\B^\ast$ is $\uw$-compact. By
Lemma~\ref{lemmaT}, $\mathbb{T}(\B^\ast)\subseteq \mathbb{T}(\B)$.
Every set of formulas $T$ from $\mathbb{T}(\B)$ is realized in
$\B^\ast$. So $\mathbb{T}(\B^\ast)=\emptyset$.
\end{proof}

\begin{corollary}
For an arbitrary algebra $\B$ there exists $\uw$-compact algebra
$\B^\ast$ which is elementary equivalent to $\B$.
\end{corollary}

In Theorem~\ref{el_uw} we constructed $\uw$-compact extension
$\B^\ast$ of $\B$ such that $\B^\ast$ is elementary equivalent to
$\B$. One can modify the idea of Theorem~\ref{el_uw} and find more
constructive $\uw$-compact extension $\B$ which is universally
equivalent to $\B$.

\begin{proposition}
Let $\B$ be an $\L$-algebra. Then there exists an extension $\Ce$
of $\B$ such that $\Ce$ is $\uw$-compact and $\Ce\equiv_\forall
\B$. Moreover, one can get $\Ce$ by (transfinite) induction in
series of extensions
$$
\B=\Ce_0 < \Ce_1 < \Ce_2 \ldots \: ,
$$
where $\Ce_{\beta+1}$ is finitely generated extension of
$\Ce_\beta$, and $\Ce_\alpha=\bigcup_{\beta < \alpha}{\Ce_\beta}$
is the union of the chain for a limit ordinal $\alpha$. Also
$\Ce_\alpha\equiv_\forall \B$ for all $\alpha$.
\end{proposition}

\begin{proof}
Let us construct $\Ce$ by means of transfinite induction on
$|\mathbb{T}(\B)|$. Take $\Ce_0=\B$. Consider an algebra $\B_1$
where $T_0$ is realized. Let $b_1,\ldots,b_n \in B_1$ be elements
such that $\B_1 \models T(b_1,\ldots,b_n)$. Put $\Ce_1$ as the
subalgebra of $\B_1$ generated by subalgebra $\B$ and elements
$b_1,\ldots,b_n$. And so on.

If $\alpha=\beta+1$ then we take $\B_\alpha$ as an ultrapower of
$\Ce_\beta$ where $T_\beta$ is realized, and $\Ce_\alpha$ is
subalgebra of $\B_\alpha$ generated by $\Ce_\beta$ and finite set
of element in $\B_\alpha$ which realize formulas from $T_\beta$.
It is easy that $\Ce_\alpha \equiv_\forall \Ce_\beta$.

For a limit ordinal $\alpha$ we put $\Ce_\alpha=\bigcup_{\beta <
\alpha}{\Ce_\beta}$ as the union of the chain $(\Ce_\beta, \beta <
\alpha)$. In this case $\Ce_\alpha=\underrightarrow{\lim} \;
\Ce_\beta$ is also the direct limit of the direct system
$(\Ce_\beta, \beta < \alpha)$, therefore $\Ce_\alpha\in
\ucl(\{\C_\beta, \beta <
\alpha\})$~\cite[Theorem~1.2.9]{Gorbunov}. Since $\Ce_\beta <
\Ce_\alpha$ we have $\Ce_\beta \in \ucl (\C_\alpha)$ for all
$\beta < \alpha$. By induction, $\Ce_\beta \equiv_\forall
\Ce_\gamma$ for any $\beta,\gamma < \alpha$. Therefore,
$\Ce_\alpha \equiv_\forall \Ce_\beta$ for every $\beta < \alpha$.

At the end of such process we get an extension $\Ce$ of $\B$ such
that $\Ce\equiv_\forall \B$ and $\mathbb{T}(\Ce)=\emptyset$, i.e.,
$\Ce$ is $\uw$-compact.
\end{proof}

The following results are also useful in universal algebraic
geometry.

\begin{lemma}\label{lemmaQ}
Let $\B,\Ce$ be an $\L$-algebras. Suppose that $\B$ is
$\qw$-compact, $\Ce\in \qvar (\B)$, and every finitely generated
subalgebra $\B_0 < \B$ is separated by $\Ce$. Then $\Ce$ is
$\qw$-compact and and $\qvar(\B)=\qvar(\Ce)$.
\end{lemma}

\begin{lemma}\label{lemmaU}
Let $\B,\Ce$ be an $\L$-algebras. Suppose that $\B$ is
$\uw$-compact, $\Ce\in \ucl (\B)$, and every finitely generated
subalgebra $\B_0 < \B$ is discriminated by $\Ce$. Then $\Ce$ is
$\uw$-compact and $\ucl(\B)=\ucl(\Ce)$.
\end{lemma}

\begin{proof}
We prove only statement about $\uw$-compactness. Statement about
$\qw$-compactness may be proven in much the same way. By
Theorem~\ref{theoremU}, it is sufficient to show that
$\ucl(\Ce)_\omega\subseteq\Dis(\Ce)_\omega$ (inclusion
$\ucl(\Ce)_\omega\supseteq \Dis(\Ce)_\omega$ holds anyway). As
$\Ce\in \ucl (\B)$ then $\ucl(\Ce) \subseteq \ucl(\B)$ and
$\ucl(\Ce)_\omega \subseteq \ucl(\B)_\omega$. Since $\B$ is
$\uw$-compact we have $\ucl(\B)_\omega=\Dis(\B)_\omega$. If every
finitely generated subalgebra $\B_0 < \B$ is discriminated by
$\Ce$ then $\Dis(\B)_\omega \subseteq \Dis(\Ce)_\omega$.
Therefore, $\ucl(\Ce)_\omega\subseteq\Dis(\Ce)_\omega$, as
desired. Also we got $\ucl(\B)_\omega=\ucl(\Ce)_\omega$ that
implies $\ucl(\B)=\ucl(\Ce)$.
\end{proof}

For $\L$-algebra $\A$ we denote by $\L_\A = \L \cup\{c_a \mid a
\in A\}$ the language $\L$ extended by elements from $\A$ as new
constant symbols~\cite[subsection~3.4]{DMR1}. An algebra $\M$ in
$\L_\A$ is called {\em $\A$-algebra} if the map $h\colon \A \to
\M$, $h(a)=c^{\:\M}_a$, $a\in A$, is embedding.

\begin{proposition}
Let $\A$ be an $\L$-algebra. Consider $\A$ as $\A$-algebra. If
$\A$ is $\qw$-compact (in the language $\L_\A$) then every
$\A$-algebra $\Ce$ from $\qvar_\A(\A)$ is $\qw$-compact. If $\A$
is $\uw$-compact (in the language $\L_\A$) then every $\A$-algebra
$\Ce$ from $\ucl_\A(\A)$ is $\uw$-compact.
\end{proposition}

\begin{proof}
We use here Lemmas~\ref{lemmaQ} and~\ref{lemmaU}. Every finitely
generated subalgebra $\A_0$ of $\A$ is an $\L_\A$-algebra, so
$\A_0=\A$. Since $\Ce$ is $\A$-algebra then $\A$ is separated and
discriminated by $\Ce$ in the obvious way. Thus, we have obtained
the looked-for result.
\end{proof}

{\small The information of the authors:

\medskip

Evelina Yu. Daniyarova, Vladimir N. Remeslennikov

Omsk Department of Institute of Mathematics, Siberian Branch of
the Russian Academy of Sciences

644099 Russia, Omsk, Pevtsova st. 13

Phone: +73812972251

e-mail: \texttt{evelina.omsk@list.ru, remesl@ofim.oscsbras.ru}

\medskip

Alexei G. Myasnikov

Department of Mathematics and Statistics, McGill University

Burnside Hall, Room 1005, 805 Shebrooke Street West, Montreal,
Quebec, Canada, H3A 2K6

Phone: +15143985476

e-mail: \texttt{amiasnikov@gmail.com}}

\end{document}